\documentclass[leqno,12pt]{article}

\usepackage{amssymb}
\usepackage{euscript}
\usepackage[dvips]{graphicx}
\usepackage{flafter}
\usepackage{pstricks}
\usepackage[french]{babel}
\usepackage[latin1]{inputenc}
\usepackage{amsmath, amssymb}
\usepackage{amsfonts}
\usepackage{pstricks-add}
\usepackage{dsfont}
\usepackage{bm}

\usepackage{mathrsfs} 

\evensidemargin\oddsidemargin

\begin{document}

\newcommand{\eqnsection}{
\renewcommand{\theequation}{\thesection.\arabic{equation}}
   \makeatletter
   \csname  @addtoreset\endcsname{equation}{section}
   \makeatother}
\eqnsection

\def\r{{\mathbb R}}
\def\e{{\mathbb E}}
\def\p{{\mathbb P}}
\def\P{{\bf P}}
\def\E{{\bf E}}
\def\Q{{\bf Q}}
\def\z{{\mathbb Z}}
\def\N{{\mathbb N}}
\def\T{{\mathbb T}}
\def\G{{\mathbb G}}
\def\L{{\mathbb L}}
\def\1{{\mathds{1}}}
\def\deg{\chi}
\def\t{{\bf{t}}}
\def\a{{\bf{a}}}
\def\C{{\mathcal{C}}}

\def\ee{\mathrm{e}}
\def\d{\, \mathrm{d}}
\def\S{\mathscr{S}}
\def\bs{{\tt bs}}


\newtheorem{theorem}{Theorem}[section]

\newtheorem{definition}[theorem]{Definition}
\newtheorem{lemme}[theorem]{Lemma}
\newtheorem{proposition}[theorem]{Proposition}

\newtheorem{corollary}[theorem]{Corollary}


\title{ First order transition for   the   branching random walk at the critical parameter}
\author{Thomas Madaule, Universit\'e Paris XIII   }
\maketitle

\vglue50pt
 
\bigskip
\bigskip

{\leftskip=2truecm \rightskip=2truecm \baselineskip=15pt \small

\noindent{\slshape\bfseries Abstract.} 
Consider a branching random walk on the real line in the boundary case. The associated additive martingales can be viewed as the partition function of a directed polymers on a disordered tree. By studying the law of the trajectory of a particle chosen under the polymer measure, we establish a first order transition for the partition function at the critical parameter. This result is strongly related to the paper of Aïdékon and Shi \cite{AShi11} in which they solved the problem of the normalisation of the partition function in the critical regime. } 

\bigskip
\bigskip

\section{Introduction}

We consider a real-valued branching random walk: Initially, a single particle sits at the origin.  Its children together with their  displacements,  form a point process $\Theta$ on $\r$ and  the first generation of the branching random walk.  These children have children of their own which form the second generation, and behave --relatively to their respective positions at birth-- like independent copies of the same point process $\Theta$.  And so on.

Let $\mathbb{T}$ be the genealogical tree of the particles in the branching random walk. Plainly, $\mathbb{T}$ is a Galton-Watson tree.  We write $|z|=n$ if a particle  $z$  is in the $n$-th generation, and denote its position by $V(z)$.  The collection of positions $(V(z),z\in \mathbb{T})$ is our branching random walk.

Following \cite{AShi11}, assume throughout the paper the following conditions
\begin{eqnarray}
\label{criticalcondition1}
\E\Big(\underset{|x|=1}{\sum}  \ee^{-V(x)} \Big)=1,&   &\E\Big(   \underset{|x|=1}{\sum}1  \Big) >1 ,\quad \text{  and }
\\
\label{criticalcondition2} \E\Big(\underset{|x|=1}{\sum}V(x)\ee^{-V(x)} \Big)=0,&   \qquad &\sigma^2:=\E \Big(\underset{|x|=1}{\sum}V(x)^2\ee^{-V(x)} \Big)<\infty.
\end{eqnarray}

The branching random walk is then said to be in the boundary case (Biggins and Kyprianou \cite{Bky05}). We refer to \cite{Jaf09} for detailed discussions on the nature of the assumption (\ref{criticalcondition1}) and (\ref{criticalcondition2})).

Let $\Phi(t):=\log \E \Big(\underset{|x|=1}{\sum}\ee^{-tV(x)} \Big)  \in (-\infty, +\infty],\, t\in \r $ and let
\begin{equation}
W_{\beta,n}:=\underset{|x|=n}{\sum}\ee^{-\beta V(x)-\Phi(\beta)n},\quad \beta\in \r,
\end{equation} 
which can be viewed as the normalized partition function of a directed polymer on trees, see the forthcoming (\ref{defmubeta}). 
In the literature $W_{1,n} $ is called the {\it critical additive martingale} associated with the branching random walk. For notational simplification, we write $W_{1,n}=W_n$ for any $n \geq 0$ [$W_0:=1$].

Under (\ref{criticalcondition1}), $\mathbb{T}$ is infinite with positive probability. Moreover the results stated here make only a trivial sense if the system dies out, so it is convenient to introduce:
\begin{eqnarray}
\P^*(\cdot):=\P\left(\cdot\, |\, \mathbb{T} \text{ is infinite}\right).
\end{eqnarray}

By  Biggins \cite{Big77},  it is known that under some  integrability conditions (for example under the forthcoming (\ref{extra1}), we refer to  Lyons \cite{Lyo97} for the optimal conditions), we have
\begin{eqnarray*}
\text{ for }   \beta<1,\qquad &&\underset{n\to\infty}{\lim}W_{\beta,n}=W_\beta>0, \qquad \P^*\text{a.s},
\\
\text{ for }  \beta \geq 1,\qquad &&\underset{n\to\infty}{\lim}W_{\beta,n}=0,\qquad\qquad \,\,\,\, \P^*\text{a.s}.
\end{eqnarray*}

According to the terminology in the study of polymers (see e.g \cite{CYos06}), we say that the region $\beta>1$ is the strong disorder regime, $\beta\leq 1$ the weak disorder regime and $\beta=1$ the critical case.

We are interested here in  the regularity of $\beta \to W_\beta$ at $\beta=1$. Biggins \cite{Big91} proved that the martingale $(W_{n,\beta})_{\beta \in \mathbb{C}}$ converges uniformly on any compact subset of a set $\Lambda^*\subset\mathbb{C}$ almost surely and in mean. As a by-product, he obtained the analyticity of $W_\beta$ on $(0,1)$.  We shall show that there is a first order transition at  $\beta=1$. In order to state our main result, we need to assume that there exist $\frac{1}{4}>\epsilon_0>0$ and $\frac{\epsilon_0}{2}>\delta_->0$ such that
\begin{eqnarray}
\label{extra1}
&& \E \Big( \Big(\underset{|x|=1}{\sum}\ee^{-(1-2\delta_-)V(x)} \Big)^{1+2\epsilon_0} \Big)<\infty.
\end{eqnarray}
Note that this condition (\ref{extra1}),     stronger than the Aïdékon-Shi \cite{AShi11}'s conditions (see (\ref{1.5})),  implies  
\begin{eqnarray}
\label{extra2}
&& \underset{\beta\in[1-\delta_-,1]}{\sup} \E\Big( \Big(\underset{|x|=1}{\sum}\ee^{-\beta V(x)} \Big)^{1+\epsilon_0} \Big)<\infty.
\end{eqnarray}
Let us also introduce the so-called {\it derivative martingale} defined by 
$$D_n:=\underset{|u|=n}{\sum}V(u)\ee^{-V(u)}, \qquad n\ge 1.$$
Defining $X:=\underset{|x|=1}{\sum}\ee^{-V(x)}$ 	and $ \tilde{X}:=\underset{|x|=1}{\sum}\max\{0,V(x)\}\ee^{-V(x)}$, Biggins and Kyprianou, \cite{BKy04}, have shown that under the condition 
\begin{eqnarray}
\label{1.5}
&&\E(X(\max(0,\log X)^2)<\infty, \qquad  \E(\tilde{X}\max (0,\log \tilde{X}))<\infty, 
\end{eqnarray}
there exists a random variable $D_\infty$ positive on the set of non-extinction such that
\begin{equation}
\label{D_infty}
\underset{n\to\infty}{\lim}D_n = D_\infty,\quad \P^* \text{a.s}.
\end{equation} 
Our first result in this paper is the following theorem:
\begin{theorem}
Assume (\ref{criticalcondition1}), (\ref{criticalcondition2}) and (\ref{extra1}). We have:
\label{mainresult}
\begin{equation}
\label{eqmainresult}
\underset{\beta\uparrow 1}{\lim}\frac{W_\beta}{1-\beta}=2D_\infty,
\end{equation}
where the convergence holds in $\P^*$ probability.
\end{theorem}
%

Theorem \ref{mainresult} relies on a study of the polymer measure at the critical point which will be our second result in this paper. Following Derrida and Spohn \cite{DSpo88}, we associate each vertex $x \in \mathbb{T}$ to $[\emptyset,x]$ the unique shortest path relating $x$ to the root $\emptyset$ and $x_i$ (for $0\leq i \leq|x|$) the vertex on $[\emptyset,x]$ such that $|x_i|=i$. The trajectory of $x\in\mathbb{T}$ corresponds to the ancestor's positions of $x$, i.e the vector $\left(V(x_1),...,V(x_{|x|})\right)$, whereas $(V_s(x))_{s\in [0,1]}$ designates the linear interpolation of the trajectory of $x\in \mathbb{T}$ and is defined by
\begin{equation}
\label{1.10} V_t(x):=\frac{1}{\sqrt{|x|}}V(x_{\lfloor |x|t\rfloor})+(|x|t-\lfloor |x|t\rfloor)\frac{1}{\sqrt{|x|}}(V(x_{\lfloor |x|t\rfloor+1})-V(x_{\lfloor |x| t \rfloor })), \quad 0\leq t \leq 1.
\end{equation}
Then for each parameter $\beta>0$ and $ n\in \N$, we define the  polymer measure  $ \mu^{(\beta)}_n $ on the disordered tree ${\mathbb T}_n$ by
\begin{equation}
\nu^{(\beta)}_n(x):= \frac{ 1}{W_{\beta, n}}\ee^{- \beta V(x) - \Phi(\beta) n}\delta_x,\qquad  x\in \mathbb{T}_n.
\end{equation}
We will study the law of the trajectory of a particle chosen under the polymer measure. So let $(\mathcal{C},||.||_\infty)$ be the set of continuous function on $[0,1]$ endowed with the sup-norm $||.||_\infty$ and for any  $A\in {\cal B}$ (the $\sigma$-algebra generated by the open sets of $(\mathcal{C},||.||_\infty)$) define  
\begin{equation}
\label{defmubeta} 
\mu^{(\beta)}_n(A):= \frac{1}{W_{\beta, n}}\underset{|x|=n}{\sum}\ee^{- \beta V(x) - \Phi(\beta) n}\1_{\{(V_s(x))_{s\in [0,1]}\in A\}}.
\end{equation}
The model of directed polymer on a disordered tree  corresponds in some sense to a mean field limit when the dimension goes to infinity, see \cite{DSpo88}. This approximation is can be understood in large lattice dimension: indeed, as $d$ increases, two independent paths $V_1$ and $V_2$ on the lattice have smaller probability to ever meet in the future. The models on the $d$-dimensional lattice and on tree with branching number $b$ are asymptotically alike  when $b=2d\to\infty$. Many authors have already worked on this subject introduced in 1988 by Derrida and Spohn \cite{DSpo88}. Recently Mörters and Ortgiese \cite{MOrt08} studied  the phase transition arising from the presence of a random disorder.  Hu and Shi in \cite{HShi09} showed that the derivative martingale appears naturally in the rate convergence of $W_n \to 0$. Furthermore, Aïdékon and Shi \cite{AShi11} proved the following theorem 
\\

\noindent {\bf Theorem A} ({\bf Aïdékon and Shi \cite{AShi11}})
{ \it Assume (\ref{criticalcondition1}), (\ref{criticalcondition2}) and (\ref{1.5})  we have
\begin{equation}
\label{eqthm1.1}
\underset{n\to\infty}{\lim} \sqrt{n} \frac{W_n}{D_n}=\left(\frac{2}{\pi\sigma^2}\right)^{\frac{1}{2}},\qquad \text{ in } \P^*\text{ probability}.
\end{equation}}
This result, used repeatedly in our paper, solved the problem of the normalisation of the partition function in the critical regime. Moreover we use the powerful method developed in \cite{AShi11} and establish our second result:
\begin{theorem}
\label{ConvL2}
Assume (\ref{criticalcondition1}), (\ref{criticalcondition2}) and (\ref{1.5}). For any  $F\in C_b(\C,\r^+)$ we have
\begin{equation}
\label{eqConvL2}
\mu_n^{(1)}(F):=\frac{1}{W_n}\underset{|u|=n}{\sum}\ee^{-V(u)}F\left((V_s(u))_{s\in [0,1]}\right) \underset{n\to\infty}{{\to}}\E\left(F(\sigma  (R_s)_{s\in[0,1]})\right),\qquad \text{in } \P^* \text{ probability,} 
\end{equation}
where $(R_s)_{s\in[0,1]}$ is a Brownian meander.
\end{theorem}
This convergence represents an important step in the proof of Theorem \ref{mainresult} but it may  also  have an independent interest. For example we mention an interesting paper by  Alberts and Ortgiese \cite{AOrt12} who also study a phase transition at the critical case. Theorem \ref{ConvL2} would yield their Theorem 1.2.

Theorem \ref{ConvL2} gives also an interesting consequence on the ``overlap" of the branching random walk which is introduced in \cite{DSpo88}. For $|u|,\,|v|=n$ we define the overlap by 
$$Q_{u,v}=\frac{\sup \left\{k\leq n,\, u_j=v_j\,\forall j\leq k\right\}}{n}.$$
Similarly for $|u|=|v|=n$ we can introduce the fraction of time in which the two paths $(V(u_1),...,V(u_n)))$ and $(V(v_1),...,V(v_n)))$ are identical, i.e
$$\tilde{Q}_{u,v}=\frac{\sup \left\{k\leq n,\, V(u_j)=V(v_j)\,\forall j\leq k\right\}}{n}.$$
Clearly, $0<Q_{u,v}\leq \tilde{Q}_{u,v}\leq 1$.

\begin{corollary}
\label{overlap}
Assume (\ref{criticalcondition1}), (\ref{criticalcondition2}) and (\ref{1.5}). For any $\delta>0$, the following convergence
 is true
\begin{equation}
\label{eqoverlap} \frac{1}{W_n^2}\underset{|u|=n,|v|=n}{\sum}\ee^{-V(u)}\ee^{-V(v)}\1_{\{ \tilde{Q}_{u,v}\geq \delta\}}\underset{t\to \infty}{\to} 0,\quad \text{in } \P^* \text{ probability.}
\end{equation}
\end{corollary}

Finally we stress that Theorem \ref{ConvL2} is not true $\P^*$ almost surely. Indeed let us introduce $F_\epsilon\in  C_b(\C,\r^+)$ defined by
\begin{eqnarray*}
F_\epsilon(w):= \left\{ \begin{array}{ll} 0& \text{ if } w(1)\notin [-\epsilon,2\epsilon],
\\
\frac{\epsilon+x}{\epsilon}& \text{ if } w(1) \in [-\epsilon,0],
\\
1& \text{ if  } w(1) \in [0,\epsilon],
\\
\frac{2\epsilon-x}{\epsilon}& \text{  if  } w(1) \in [\epsilon,2\epsilon],
\end{array}\right. \qquad \epsilon>0,\, w\in \mathcal{C},
\end{eqnarray*}
and
\begin{eqnarray*}
A_n:= \{\exists x:\,n\leq |x|\leq 2n,\, \frac{1}{2}\log n \leq V(x)\leq \frac{1}{2}\log n+C;\,\underset{n\leq k \leq2n}{\max}\sqrt{n} W_k\leq C \},\quad n\in \N,\, C>0.
\end{eqnarray*}
According to Lemma 6.3 in \cite{AShi11} and Theorem 1.5 in \cite{HShi09}, there exist $ C,\, N>0$ large and $c>0$ small, such that for any $n>N$, $\P^*(A_n)\geq c$. Clearly for any $\epsilon>0,\, n\geq \epsilon^{-2}$, on the set $A_n$, we have $ \mu_n^{(1)}(F_\epsilon) \geq C^{-1}\ee^{-C}.$ So
\begin{eqnarray*}
\underset{\epsilon\to 0}{\limsup}\, \underset{n\to\infty}{\limsup}\, \E^*(\mu_n{(1)}(F_\epsilon))\geq cC^{-1}\ee^{-C} > 0=\underset{\epsilon\to 0}{\limsup}\, \E\left( F_\epsilon[(\sigma R_s)_{s\in [0,1]}]\right),
\end{eqnarray*}
which implies that (\ref{eqConvL2}) can not hold $\P^*$ almost surely.
\\

Similarly, (\ref{eqoverlap}) can not be strengthened in $\P^*$ almost sure convergence.
\\

The rest of the paper is organized as follows. In Section 2 we present some  preliminaries on branching random walks. Section 3 is devoted to the proof of the Theorem \ref{ConvL2} and Corollary \ref{overlap}. Finally,  in Section 4 we prove the Theorem \ref{mainresult}.

\section{Preliminaries}

This section   collects    some preliminary results on the branching random walk (change of probabilities, an associated one-dimensional random walk),  and it entirely  comes from   Aïdékon and Shi  \cite{AShi11}.

\subsection{The many-to-one Lemma}
Let $(V(x))$ be a branching random walk satisfying (\ref{criticalcondition1}) and (\ref{criticalcondition2}). Let $(S_n)_{n\geq 0}$ a random walk such that the law of $(S_1)$ is given by 
\begin{equation}
\label{defSn}
\E\left(f(S_1)\right):=\E \Big(\underset{|z|=1}{\sum}f(V(z))\ee^{-V(z)} \Big),\qquad \forall f:\r\to [0,\infty),\,\,\text{measurable}.
\end{equation}
The condition (\ref{criticalcondition1}) and (\ref{criticalcondition2}) implies that $(S_n)$ is a mean zero random walk and $\E(S_1^2)=\sigma^2<\infty$. From a simple induction it stems that for any $n\geq 0$ and $g:\r^n\to \r^+$ measurable we have:
\begin{equation}
\label{2.1}
\E \Big(\underset{|x|=n}{\sum}g\left(V(x_1),...,V(x_n)\right) \Big)=\E\left(\ee^{S_n}g(S_1,...,S_n)\right).
\end{equation}
Equality (\ref{2.1}) forms the so called many-to-one Lemma which plays a fundamental role in many computations of expectations. The presence of the random walk $(S_i)$ is explained in Lyons, Pemantle, Peres \cite{LPP95}, Lyons \cite{Lyo97} and Biggins, Kyprianou \cite{BKy04}.

\subsection{The renewal function associated with a one-dimensional random walk}
Associated to $(S_n)$, which is a centered random walk real-valued with $\sigma^2=\E[S_1^2]\in(0,\infty)$, let $h_0$ be its renewal function defined by
\begin{equation}
h_0(u):=\underset{j\geq 0}{\sum}\P\left( \underset{i\leq j-1}{\min}S_i>S_j\geq -u  \right),\qquad u\geq 0.
\end{equation}
In the following this function will play an important role, so we collect here some facts about $h_0$. For any $u\geq 0$, $h_0$ satisfies
\begin{equation}
\label{3.40}
 h_0(u)=\E\left(h_0(S_1+u)\1_{\{S_1\geq -u\}}\right).
\end{equation}
If we write 
\begin{equation}
\underline{S}_n:=\underset{j\leq n}{\min}\,\, S_j,\qquad   n\geq 0,
\end{equation} 
it is known that there exists $c_0>0$ and $\theta>0$ such that
\begin{equation}
\label{lldda}
 c_0:=\underset{u\to\infty}{\lim}\frac{h_0(u)}{u},\qquad \P\left(\underline{S}_n\geq -u\right)\underset{n\to\infty}{\sim}\frac{\theta h_0(u)}{n^\frac{1}{2}},\qquad \forall u\geq 0.
\end{equation}
As a consequence there exists constants $c,\,C>0$ such that  
\begin{equation}
\label{2.7}
c_1(1+u)\leq h_0(u)\leq C_1(1+u),\qquad u\geq 0.
\end{equation}
As in \cite{AShi11}, we will need the following uniform version of (\ref{lldda}):  as $n\to \infty$,
\begin{equation}
\label{lldda2}
\P\left(\underline{S}_n\geq -u\right)= \frac{\theta h_0(u)+ o(1)}{n^\frac{1}{2}},
\end{equation}
uniformly in $u\in [0, (\log n)^{30}]$. 

Finally we mention the inequality due to \cite{AShi11}: {\it there exists $c>0$ such that for $u>0,a\geq0,\,b\geq 0$ and $n\geq 1$,
\begin{equation}
\label{Aine}
\P\left(\underline{S}_n\geq -a,\,b-a\leq S_n\leq b-a+u\right)\leq c\frac{(u+1)(a+1)(b+u+1)}{n^{\frac{3}{2}}}.
\end{equation}}


%

\subsection{A spine conditioned to stay positive}
Let $(V(x))$ be a branching random walk satisfying (\ref{criticalcondition1}) and (\ref{criticalcondition2}), let $(\mathcal{F}_n)$ be the sigma-algebra generated by the branching random walk in the first $n$ generations. Since Lyons \cite{Lyo97}, the spinal decomposition is a widespread technique to study the branching random walk. Usually we introduce the martingale $W_n:=\underset{|z|=n}{\sum}\ee^{-V(z)}$ to define the probability $\Q$ satisfying for any $n\in \N$, $\Q|_{\mathcal{F}_n}:=W_n.\P|_{\mathcal{F}_n}$. Then we obtain a description with a spine of the  branching random walk under $\Q$ , moreover this spine behaves like a centered random walk. Here we will need a slightly different decomposition, we will work with a spine whose the law is as a random walk conditioned to stay positive.

 First let us introduce some notations. For any vertex $x\in \T$, let $\underline{V}(x):=\underset{y\in]\emptyset,x]}{\min}V(y)$. Then for $\alpha\geq0$ and $u\geq -\alpha$ let $h_\alpha(u):=h_0(u+\alpha).$
Finally we define the processes
\begin{eqnarray*}
&&W_n^{(\alpha)}:=\underset{|x|=n}{\sum}\ee^{-V(x)}\1_{\{\underline{V}(x)\geq -\alpha\}}, \quad D_n^{(\alpha)}:=\underset{|x|=n}{\sum}h_\alpha(V(x))\ee^{-V(x)}\1_{\{\underline{V}(x)\geq -\alpha\}}.
\end{eqnarray*}
 From (\ref{3.40}) and the branching property stem that for any $\alpha\geq 0$, $(D_n^{(\alpha)},n\geq 0)$ is a non-negative martingale with respect to  $\mathcal{F}_n$ (see Biggins Kyprianou \cite{BKy04} or \cite{AShi11} for a proof). So associated with $D_n^{(\alpha)}$ we introduce the new probability measure $\Q^{(\alpha)}$ which satisfies for any $n$,
\begin{equation}
\label{ChangQ}\Q^{(\alpha)}|_{\mathcal{F}_n}:=\frac{D_n^{(\alpha)}}{h_\alpha(0)}.\P|_{\mathcal{F}_n}.
\end{equation}
Now we will give a representation with spine of the branching random walk under $\Q^{(\alpha)}$. A justification of this representation can be founded in \cite{AShi11}. Recall that the point process which governs the law at the first generation of $(V(x),|x|=1)$ is distributed under $\P$ as the point process $\Theta$. For any $u\geq -\alpha$, the new probability $\Q^{(\alpha)}$ makes appear the point process $\hat{\Theta}_u^{(\alpha)}$ whose distribution is the law of $(u+V(w),|x|=1)$ under $\Q^{(u+\alpha)}$. Then the branching random walk under $\Q^{(\alpha)}$ is governed by the followings rules:

\begin{itemize}
\item $w_0^{(\alpha)}=\emptyset$ gives birth to particles distributed according to $\Theta_{0}^{(\alpha)}$.
\item Choose $w_1^{(\alpha)}$ among children of $w_0^{(\alpha)}$ with probability proportional to $\ee^{-V(x)}\1_{\{V(x)\geq -\alpha\}}h_\alpha(V(x))$.
\item $\forall n\geq 1$, $w_n^{(\alpha)}$ gives birth to particles distributed according to $\Theta_{u}^{(\alpha)}$ ($u=V(w_n^{(\alpha)})$).
\item Choose $w_{n+1}^{(\alpha)}$ among the children of $w_n^{(\alpha)}$ with probability proportional to $\ee^{-V(y)}\1_{\{V(y)\geq -\alpha \}} h_{\alpha}(V(y))$.
\item Subtrees rooted at all other brother particles are independent branching random walks under $\P$.
\end{itemize}
\paragraph
\noindent See below three facts which we will use continuously:

{\it
(i) $\Q^{(\alpha)}(non-extinction)=1$ and $\forall n\in\N$, $\Q^{(\alpha)}(D_n^{(\alpha)}>0)=1$.

(ii) For any $n$ and any vertex $x$ with $|x|=n$, we have
\begin{equation}
\label{eqepine}
\Q^{(\alpha)}(w_n^{(\alpha)}=x|\mathcal{F}_n)=\frac{h_\alpha(V(x))\ee^{-V(x)}\1_{\{\underline{V}(x)\geq -\alpha\}}}{D_n^{(\alpha)}}.
\end{equation}

(iii) The spine process $(V(w_n^{(\alpha)}), n\geq 0)$ under $\Q^{(\alpha)}$, is distributed as a Markov chain with transition probabilities given by
\begin{equation}
p^{(\alpha)}(u,dv):=\1_{\{v\geq -\alpha\}}\frac{h_\alpha(v)}{h_\alpha(u)}p(u,dv),\qquad u\geq -\alpha,
\end{equation}
where $p(u,dv):=\P(S_1+u\in dv)$ is the transition probability of $(S_n)$. In the sense of Doob's $h$-transform, $(V(w_n^{(\alpha)}))_{n\in \N} $ has the law of the random walk $(S_n)_{n\geq 0}$ conditioned to stay in $[-\alpha,\infty]$. A convey way to represent this processes is the following identity: for any $n\geq 1$ and any measurable function $g:\r^{n+1}\to [0,\infty)$,
\begin{equation}
 \label{eq3.3}
\E_{\Q^{(\alpha)}}\left(g(V(w_i^{(\alpha)}),\, 0\leq i\leq n)\right)=\frac{1}{h_\alpha(0)}\E\left(g(S_i,\,0\leq i\leq n)h_\alpha(S_n)\1_{\{\underline{S}_n\geq -\alpha\}}\right).
\end{equation}}

{\bf Convention:} Throughout the paper, $c,\, c',\, c''$ denote generic constants which may change from paragraph to paragraph, but are independent of $n$.

\section{Proof of Theorem \ref{ConvL2} and Corollary \ref{overlap}}

Let us introduce $(R_s)_{s\in [0,1]}$ and $(R'_s)_{s\in [0,1]}$ two independent Brownian meander under $\P$. For any $d\in \N^*$, ${\bf t}=(t_1,...,t_d)\in [0,1]$ and any process $(f_s)_{s\in [0,1]}$ we will denote the vector $(f_{t_1},...,f_{t_d})$ by $f_{\bf t}$.

This section is divided in two steps:

- A) We show Theorem \ref{ConvL2} and Corollary \ref{overlap} assuming the following assertion: {\it Under the integrability conditions (\ref{criticalcondition1}), \ref{criticalcondition2}) and (\ref{1.5}), for any $d\in \N^*$, ${\bf t}=(t_1,...,t_d)\in[0,1]^d$, $F\in \C_b(\r^d)$,
\begin{equation}
\label{1particule}
\underset{n\to\infty}{\lim} \frac{1}{W_n}\underset{|u|=n}{\sum}\ee^{-V(u)}   F(V_{\bf t}(u))=\E\left(F(\sigma R_{{\bf t}})\right),\quad \text{in } \P^* \text{ probability.}
\end{equation}}

- B) We prove assertion (\ref{1particule}).

\subsection{ Step A)}
For any metric spaces $E$ and $F$ we denote $\mathcal{C}_b(E,F):=\{f: E\to F,\, \text{continuous and bounded}\}$.  Let $\mathbb{W}^+$ the law of the Brownian meander, and $(\E^*(\mu_n\otimes\mu_n(\cdot)))_{n\geq 0}$ the sequence of probability measure on $\mathcal{C}^2:=\{f:[0,1]\to \r^2,\, \text{continuous}\}$, defined by
\begin{equation}
\label{defmeasureprod}
  \E^*(\mu_n\otimes\mu_n(F)):= \E^*\Big(\frac{1}{W_n^2} \sum_{|u|=|v|=n}\ee^{-V(u)}\ee^{-V(v)}F\left[ (V_t(u),V_t(v))_{t\in [0,1]}\right]\Big),\qquad  F\in  C_b(\C^2,\r).
\end{equation}
First we shall prove that (\ref{1particule}) implies:
\begin{equation}
\label{Weak}
\E^*(\mu_n\otimes\mu_n(\cdot)) \overset{weakly}{\Longrightarrow} \mathbb{W}^+\otimes \mathbb{W}^+.
\end{equation}
In order to obtain (\ref{Weak}), recall from \cite{Bil99} that for any continuous process the convergence of the finite-dimensional laws and the relative compactness imply the weak convergence. We will obtain the relative compactness via the following criteria (see \cite{RYor99}): 

{\it A sequence $(\P_n)$ of probability measures on $\C^2$ is weakly relatively compact if and only if the following two conditions hold:

i) for every $\epsilon>0$, there eixst a number $A$ and an integer $n_0$ such that
$$\P_n[|w(0)|>A]\leq \epsilon,\qquad \text{ for every }n\geq n_0.$$

ii) for every $\eta,\,\epsilon>0$, there exists a number $\delta$ and an integer $n_0$
$$\P_n[K(.,\delta)>\eta]\leq \epsilon,\qquad \text{ for every }n\geq n_0.$$
with
$$K(w,\delta):=\sup\left\{|w_1(t)-w_1(t')|+|w_2(t)-w_2(t')|;\, |t-t'|\leq \delta \right\},\qquad \forall \delta>0,\, w=(w_1,w_2)\in \C^2.$$ }

{\it Proof of the relative compactness of $(\E^*(\mu_n\otimes\mu_n(\cdot)))$.}
The first condition is trivially satisfied. For the second we need to control $\E^*\big(\frac{1}{W_n^2}\underset{|u|=n,|v|=n}{\sum}\ee^{-V(u)-V(v)}\1_{\{K(V(u,.),V(v,.))\geq \eta\}}\big)$. It is known that $\underset{A\to\infty}{\lim}\underset{n\in\N}{\sup}\,\, \P^*\left(\frac{1}{W_n}\geq A\sqrt{n}\right)= \underset{\alpha\to\infty}{\lim} \P^*\big(\underset{x\in\T}{\min}V(x)\leq -\alpha \big)=0$, see for instance (\ref{eqthm1.1}) and (1.4) in \cite{AShi10}. Thus for any $\epsilon>0$ there exist $A,\,\alpha>0$ large enough such that for any  $n\in \N^*$, we have:

\begin{eqnarray}
\label{Rcompact}\E^*_{(\ref{Rcompact})}&:=&\E^*\Big(\frac{1}{W_n^2}\underset{|u|=n,|v|=n}{\sum}\ee^{-V(u)-V(v)}\1_{\{K(V_{\cdot}(u),V_\cdot(v))\geq \eta\}}\Big)
\\
\nonumber &\leq&  \E^*\Big(\frac{1}{W_n}\underset{|u|=n}{\sum}\ee^{-V(u)}\1_{\{\sup\left\{|V_t(u)-V_{t'}(u)|;\, |t-t'|\leq \delta \right\}\geq \frac{\eta}{2}\}}\Big)
\\
\nonumber &\leq&  \epsilon + cA\sqrt{n}\E\Big(\underset{|u|=n}{\sum}\ee^{-V(u)}\1_{\{\underline{V}(u)\geq -\alpha\}}\1_{\{\sup\{|V_t(u)-V_{t'}(u)|;\, |t-t'|\leq \delta   \}\geq \frac{\eta}{2}\}}\Big).
\end{eqnarray}
Now using the Many-to-one Lemma  we have
\begin{eqnarray*}
\E^*_{(\ref{Rcompact})} &\leq &  \epsilon + cA\sqrt{n} \P\Big( \sup\left\{|S(n,t)-S(n,t')|;\,|t-t'|\leq \delta \right\}\geq \frac{\eta}{2};\underline{S}_n\geq -\alpha\Big),
\end{eqnarray*}
with $S(n,t):= \frac{1}{\sqrt{n}} S_{\lfloor nt \rfloor}+ (nt-\lfloor nt \rfloor)(S_{\lfloor nt \rfloor +1} -S_{\lfloor nt \rfloor}),\, t\in [0,1]$. Finally as  
\begin{eqnarray*}
\underset{\delta\to 0}{\lim}\,\underset{n\to\infty}{\limsup}\, \sqrt{n}\P\left(\sup\left\{|S(n,t)-S(n,t')|;\,|t-t'|\leq \delta \right\}\geq \frac{\eta}{2};\underline{S}_n\geq -\alpha \right)=0,\quad \text{(cf \cite{Igl74} pp 615)},
\end{eqnarray*}
we deduce that for any $\eta,\, \epsilon>0$ there exist $\delta>0,\, N\in \N$ such that for any $n\geq N$, $\E^*_{(\ref{Rcompact})}\leq 3\epsilon$. It ends the proof of the relative compactness of $(\E^*(\mu_n\otimes\mu_n(\cdot)))$. \hfill$\Box$
\\

Now we have to prove that the finite dimensional distributions of $\E^*(\mu_n\otimes \mu_n(\cdot))$ converge to those of $\mathbb{W}^+\otimes\mathbb{W}^+$, i.e:  {\it for any $d\in \N^*$, ${\bf t}=(t_1,...,t_d)\in [0,1]^d$ and $F\in C_b(\r^{2d},\r)$,
\begin{equation}
\label{finidim}
 \underset{n\to\infty}{\lim}\E^*\Big(\frac{1}{W_n^2} \sum_{|u|=|v|=n} \ee^{-V(u)-V(v)}F(V_{\bf t}(u),V_{\bf t}(v))\Big)= \E\Big(F(\sigma R_{{\bf t}},\sigma R'_{\bf t})\Big).
\end{equation} }
\\

According to (\ref{1particule}), for any ${\bf t}=(t_1,...,t_d)\in [0,1]^d$ and $(F,G)\in (C_b(\r^{d},\r))^2$
\begin{eqnarray}
\label{tou0}
\underset{n\to\infty}{\lim} \frac{1}{W_n}\underset{|u|=n}{\sum}\ee^{-V(u)} F(V_{\bf t}(u)) =\E\left(F(\sigma R_{{\bf t}})\right)  ,\quad \text{in } \P^* \text{ probability,}
\\
\label{tou1}
\underset{n\to\infty}{\lim} \frac{1}{W_n}\underset{|v|=n}{\sum}\ee^{-V(v)}   G(V_{\bf t}(v)) =\E\left(G(\sigma R_{{\bf t}})\right),\quad \text{in } \P^* \text{ probability.}
\end{eqnarray}
The left hand terms of (\ref{tou0}) and (\ref{tou1}) are bounded, then by taking the expectation of the product of (\ref{tou0}) and (\ref{tou1}) we get
\begin{eqnarray*}
\underset{n\to\infty}{\lim}\E^*\Big( \frac{1}{W^2_n}\underset{|v|=|u|=n}{\sum}\ee^{-V(u)-V(v)}   F(V(u,{\bf t}))G(V(v,{\bf t}))\Big )=\E\left(F(\sigma R_{{\bf t}})\right) \E\left(G(\sigma R'_{{\bf t}})\right).
\end{eqnarray*}
This equality is sufficient to affirm that the law of $ (V(u,{\bf t}),(V(v,{\bf t}))$ under $\E^*(\mu_n\otimes \mu_n(\cdot))$ converges to this one $(R_{\bf t},R_{\bf t}')$ under $\P$, which implies (\ref{finidim}).
\hfill$\Box$
\\

To conclude {\bf step A)} it remains to show that (\ref{Weak}) implies Theorem \ref{ConvL2} and Corollary \ref{overlap}. For any $F\in C_b(\mathcal{C},\r)$ let $\mathtt{F}_{2*} \in C_b(\mathcal{C}^2,\r)$ be the function defined by: 
\begin{eqnarray*}
&&\mathtt{F}_{2*}  ( (w_1(t),w_2(t))_{t\in [0,1]}):=\left[F((w_1(t))_{t\in [0,1]})-\E(F((R_s)_{s\in [0,1]})) \right]\times
\\
&&\qquad \qquad\qquad\qquad   \qquad\qquad \qquad\left[F(( w_2(t))_{t\in [0,1]})-\E(F((R_s)_{s\in [0,1]})) \right],\qquad    w_1,w_2\in \mathcal{C}.
\end{eqnarray*}
Then according to (\ref{Weak}), for any $F\in C_b(\mathcal{C},\r)$ we have
\begin{eqnarray*}
 \E^*\left(\left(\mu_n(F)-\E[F(R)]\right)^2\right)=\E^*(\mu_n\otimes\mu_n(\mathtt{F}_{2*}))\to \mathbb{W}^+\otimes\mathbb{W}^+(\mathtt{F}_{2*})=0,
\end{eqnarray*}
which implies Theorem \ref{ConvL2}. \hfill$\Box$
\\

Concerning Corollary \ref{overlap}, observing that $\{\tilde{Q}_{u,v}\geq \delta\}\subset \{V(u,s)=V(v,s),\,\forall s\leq \delta\}$, we deduce that (\ref{Weak}) implies
 $$\underset{n\to\infty}{\lim}\sup\E^*\Big(\frac{1}{W_n^2}\underset{|u|=n,|v|=n}{\sum}\ee^{-V(u)}\ee^{-V(v)}\1_{\{ \tilde{Q}_{u,v}\geq \delta\}}\Big)\leq \P \left(R_s=R'_s,\,\forall s\leq \delta/2\right)=0,$$
which gives Corollary \ref{overlap}.\hfill$\Box$
\\

So we can turn now to the proof of (\ref{1particule}).

\subsection{Step B): proof of (\ref{1particule})}
Fix $d\in \N^*$.  Recall that for any $(t_1,...,t_d)\in [0,1]^d$ we  denote 
\begin{eqnarray*}
{\bf t}:= (t_1,...,t_d),\quad R_{\bf t}:= (R_{t_1},...,R_{t_d}),\quad\text{and }\,\,   V_{\bf t}(u):=(V_{t_1}(u),...,V_{t_k}(u)).
\end{eqnarray*}
Let us also introduce for any $\alpha>0,\, F\in \mathcal{C}_{b}(\r^d,\r),\,{\bf t}\in [0,1]^d,\, n\in \N$ and ${\bf y}\in \r^d$, 
\begin{equation}
\overline{F}_{\bf t} ({\bf y}):=F({\bf y})-\E(F(\sigma R_{\t})), \text{  and   }\,\,    W_n^{(\alpha),\overline{F}_{\bf t}}:=\underset{|x|=n}{\sum}\ee^{-V(x)}\1_{\{\underline{V}(x)\geq -\alpha\}}\overline{F}_{\bf t}(V_{\bf t}(x)).
\end{equation}

 Following \cite{AShi11}, to prove (\ref{1particule}), we firstly show the following result on $\Q^{(\alpha)}$: {\it for any $\alpha\geq 0$, ${\bf t}\in \r^d_+$ and $F\in \C_b(\r^d,\r)$,
\begin{equation}
\label{L2sousQ}
\underset{n\to\infty}{\lim} \E_{\Q^{(\alpha)}}\Big(\frac{1}{{(W^{(\alpha)}_n)}^2} \Big(\underset{|u|=n}{\sum}\ee^{-V(u)}\1_{\{ \underline{V}(u)\geq -\alpha\}}  \overline{F}_{\bf t}(V_{\bf t}(u))\Big)^2 \Big)=0.
\end{equation}}
Equality (\ref{L2sousQ}) represents the exact analogue of (\ref{1particule}) under $\Q^{(\alpha)}$. Because of the relation (\ref{ChangQ}) we will see at the end of this section how to obtain (\ref{1particule}) from (\ref{L2sousQ}) by letting $\alpha$ goes to infinity. Working under $\Q^{(\alpha)}$ presents the following advantage: under $\Q^{(\alpha)}$ the spine remains above a barrier positioned at $-\alpha$, then the random variables $(W_n^{(\alpha)},D_n^{(\alpha)})$ are much more concentrated around their mean than $W_n$ and  $D_n$ under $\P$.

Let us start the 

\noindent{\it Proof of (\ref{L2sousQ}).} As in \cite{AShi11} we need to rewrite $W_n^{(\alpha)},D_n^{(\alpha)}$ and $W_n^{(\alpha),\overline{F}_{\bf t}}$ according to the position of spine $V(w_n^{(\alpha)})$. For each vertex $x$ with $|x|=n$ and $x\neq w_n^{(\alpha)}$, there is a unique $i$ with $0\leq i<n$ such that $w_i^{(\alpha)}\leq x$ and that $w_{i+1}^{(\alpha)}\nleqslant x$. For any $i\geq 1$, let
$$R_i^{(\alpha)}:=\left\{|x|=i: x>w_{i-1}^{(\alpha)},x\neq w_i^{(\alpha)}\right\},$$
(in words, $R_i^{(\alpha)}$  stands for the set of "brothers'' of $w_i^{(\alpha)}$). Accordingly,
$$W_n^{(\alpha),\overline{F}_{\bf t}}=\ee^{-V(w_n^{(\alpha)})}\overline{F}_{\bf t}(V_{\bf t}(w_n^{(\alpha)}))+\overset{n-1}{\underset{i=0}{\sum}}\underset{y\in R_{i+1}^{(\alpha)}}{\sum}\underset{|x|=n,x\geq y}{\sum}\ee^{-V(x)}\1_{\{\underline{V}(x)\geq -\alpha\}}\overline{F}_{\bf t}(V_{\bf t}(x)).$$
Let $k_n<n$ be an integer such that $k_n\to \infty$ ($n\to\infty$). We write
\begin{eqnarray*}
&&W_n^{(\alpha), \overline{F}_{\bf t},[0,k_n)}:=\overset{k_n-1}{\underset{i=0}{\sum}}\underset{y\in R_{i+1}^{(\alpha)}}{\sum}\underset{|x|=n,x\geq y}{\sum}\ee^{-V(x)}\1_{\{\underline{V}(x)\geq -\alpha\}}\overline{F}_{\bf t}(V_{\bf t}(x)),
\\
&&W_n^{(\alpha), \overline{F}_{\bf t},[k_n,n]}:=\ee^{-V(w_n^{(\alpha)})}\overline{F}_{\bf t}(V_{\bf t}(w_n^{(\alpha)} ))+\overset{n-1}{\underset{i=k_n}{\sum}}\underset{y\in R_{i+1}^{(\alpha)}}{\sum}\underset{|x|=n,x\geq y}{\sum}\ee^{-V(x)}\1_{\{\underline{V}(x)\geq -\alpha\}}\overline{F}_{\bf t}(V_{\bf t}(x)),
\end{eqnarray*}
so that 
\begin{eqnarray}
\label{decompos1} W_n^{(\alpha), \overline{F}_{\bf t}}&=&W_n^{(\alpha), \overline{F}_{\bf t},[0,k_n)}+W_n^{(\alpha), \overline{F}_{\bf t},[k_n,n]}, \qquad \text{and similarly we can write},
\\
\label{decompos2} W_n^{(\alpha)}&=& W_n^{(\alpha),[0,k_n)}+  W_n^{(\alpha),[k_n,n]},\qquad  D_n^{(\alpha)}= D_n^{(\alpha),[0,k_n)}+  D_n^{(\alpha),[k_n,n]}.
\end{eqnarray}
 Let also
\begin{eqnarray}
\label{4.31}&&E_{n,1}:=\{k_n^{1/3}\leq V(w_{k_n}^{(\alpha)})\leq k_n\}\cap \overset{n}{\underset{i=k_n}{\bigcap}}\{V(w_i^{(\alpha)})\geq k_n^{\frac{1}{6}}\},
\\
\label{4.32} &&E_{n,2}:=\overset{n-1}{\underset{i=k_n}{\bigcap}}\Big\{\underset{y\in R^{(\alpha)}_{i+1}}{\sum}[1+(V(y)-V(w_i^{(\alpha)}))^+]\ee^{-[V(y)-V(w_i^{\alpha})]}\leq \ee^{V(w_i^{(\alpha)})/2}\Big\},
\\
\label{4.3} &&E_{n,3}:=\left\{D_n^{(\alpha),[k_n,n]}\leq \frac{1}{n^2}\right\},\quad \text{and  }\,\, E_n:=E_{n,1}\cap E_{n,2}\cap E_{n,3}.
\end{eqnarray}

Under $\Q^{(\alpha)}$, $(V(w_n^{(\alpha)}))_{n\geq0}$ has the law of a centered random walk conditioned to stay positive. Moreover it is well known that a such process ``tends to infinity" when $n$ goes to infinity. Then keeping this fact in minds we are inclined to affirm that $W_n^{(\alpha),\overline{F}_{\bf t},[k_n,n]}$ $D_n^{(\alpha),[k_n,n]}$ have a negligible contribution in (\ref{decompos1}) and (\ref{decompos2}). The following Lemma makes rigorous this affirmation,
\begin{lemme}[Aïdékon and Shi \cite{AShi11}]
\label{9.1}
Let $\alpha\geq 0$. Let $k_n$ be such that $\frac{k_n}{\log^6 n}\to \infty$ and that $\frac{k_n}{n^{\frac{1}{2}}}\to 0$, $n\to \infty$. Let $E_n$ be as in (\ref{4.3}). Then 
\begin{equation}
\label{eq9.1}
\underset{n\to \infty}{\lim}\Q^{(\alpha)}(E_n)=1,\qquad \underset{n\to \infty}{\lim}\underset{u\in [k_n^{\frac{1}{3}},k_n]}{\inf} \Q^{(\alpha)}\left(E_n|V(w_{k_n}^{(\alpha)})=u\right)=1.
\end{equation}
\end{lemme}
Let $(k_n)_{n\in \N^*}:= (\lfloor (\log n)^{30} \rfloor)_{n\in \N^*}$. Fix $\epsilon,\, \alpha>0$, ${\bf t}\in [0,1]^d$ and $F\in C_b(\r^d,\r)$. As $||F||<\infty$, for any event $\mathtt{A}$, the expectation in (\ref{L2sousQ}) (:=$\E_{\Q^{(\alpha)}}^{(\ref{L2sousQ})}$) satisfies
\begin{eqnarray}
\label{righthan} \E_{\Q^{(\alpha)}}^{(\ref{L2sousQ})}\leq ||F|| \Q^{(\alpha)}(\mathtt{A})+\E_{\Q^{(\alpha)}}\Big(\frac{1}{({W^{(\alpha)}_n})^2}\big(\underset{|u|=n}{\sum}\ee^{-V(u)} \1_{\{\underline{V}(u)\geq -\alpha\}} \overline{F}_{\bf t}(V_{\bf t}(u)) \big)^2\1_{\mathtt{A}}\Big)  .
\end{eqnarray}
For any $n\in \N$ we set $ A_n:=\left\{ \left|\frac{D_n^{(\alpha)}}{W_n^{(\alpha)}}-\frac{\sqrt{n}}{\theta} \right| < \sqrt{n} \right\} $. According to Proposition 4.1 of \cite{AShi11} we have $ \underset{n\to\infty}{\lim}\Q^{(\alpha)}(A_n)=1$. Moreover on $A_n$, $\frac{1}{(W_n^{(\alpha)})^2}\leq n\frac{(1+1/\theta)^2}{(D_n^{(\alpha)})^2}$, thus the expectation in the right hand term of (\ref{righthan}) is smaller than:
\begin{eqnarray}
\label{boubou}&& cn\E_{\Q^{(\alpha)}}\Big(\frac{W_n^{(\alpha),\overline{F}_{\bf t}}}{D_n^{(\alpha)}}   \1_{A_n} \underset{|u|=n}{\sum}\ee^{-V(u)} \1_{\{\underline{V}(u)\geq -\alpha\}}\frac{h_\alpha(V(u))}{D_n^{(\alpha)}}  \frac{\overline{F}_{\bf t}(V_{\bf t}(u))}{h_\alpha(V(u))} \Big)
\\
\nonumber &&= cn\E_{\Q^{(\alpha)}}\Big(\1_{A_n}  \frac{\overline{F}_{\bf t}(V_{\bf t}(w_n^{(\alpha)} ))}{h_\alpha(V(w_n^{(\alpha)}))}\frac{W_n^{(\alpha),\overline{F}_{\bf t}}}{D_n^{(\alpha)}}\Big),
\end{eqnarray}
where we have recognize in (\ref{boubou}) the expression of  $\Q^{(\alpha)}\left(w_n^{(\alpha)}=u |\mathcal{F}_n\right)$ as described in (\ref{eqepine}). Finally there exists $n_0>0$ such that for any $n\geq n_0$,
\begin{equation}
\label{tohtho}
\E_{\Q^{(\alpha)}}^{(\ref{L2sousQ})}\leq  \epsilon +  cn\E_{\Q^{(\alpha)}}\Big(\1_{A_n}  \frac{\overline{F}_{\bf t}(V_{\bf t}(w_n^{(\alpha)} ))}{h_\alpha(V(w_n^{(\alpha)}))}\frac{W_n^{(\alpha),\overline{F}_{\bf t}}}{D_n^{(\alpha)}}\Big).
\end{equation}

Now by using Lemma \ref{9.1} we shall prove that we can replace $\frac{W_n^{(\alpha),\overline{F}_{\bf t}}}{D_n^{(\alpha)}}$ by $\frac{W_n^{(\alpha),\overline{F}_{\bf t},[0,k_n)}}{D_n^{(\alpha),[0,k_n)}}$ in (\ref{tohtho}). We define
\begin{eqnarray}
\label{Diane1} I_n &:=&n\E_{\Q^{(\alpha)}}\Big(\frac{W_n^{(\alpha),\overline{F}_{\bf t}}}{D_n^{(\alpha)}}\1_{A_n}  \frac{\overline{F}_{\bf t}(V_{\bf t}(w_n^{(\alpha)} ))}{h_\alpha(V(w_n^{(\alpha)}))}\Big),\qquad \text{ and}
\\
\label{Diane2}  L_n &:=& n\E_{\Q^{(\alpha)}}\Big(\frac{W_n^{(\alpha),\overline{F}_{\bf t},[0,k_n)}}{D_n^{(\alpha),[0,k_n)}}\1_{\{V(w_{k_n}^{(\alpha)})\in [k_n^{\frac{1}{3}},k_n]\}}\frac{\overline{F}_{\bf t}(V_{\bf t}(w_n^{(\alpha)} ))} {h_\alpha(V(w_n^{(\alpha)}))}\Big).
\end{eqnarray}
In the following we shall prove that $|I_n-L_n|\underset{n\to \infty}{\to} 0$. To achieve this goal we also introduce
\begin{eqnarray}
\label{Diane3}  J_n &:=&n\E_{\Q^{(\alpha)}}\Big(\frac{W_n^{(\alpha),\overline{F}_{\bf t}}}{D_n^{(\alpha)}}\1_{E_n}\frac{\overline{F}_{\bf t}(V_{\bf t}(w_n^{(\alpha)} ))} {h_\alpha(V(w_n^{(\alpha)}))}\Big),\qquad \text{and}
\\
\label{Diane4}  K_n &:=& n\E_{\Q^{(\alpha)}}\Big(\frac{W_n^{(\alpha),\overline{F}_{\bf t},[0,k_n)}}{D_n^{(\alpha),[0,k_n)}}\1_{E_n}\frac{\overline{F}_{\bf t}(V_{\bf t}(w_n^{(\alpha)} ))}{h_\alpha(V(w_n^{(\alpha)}))}\Big),
\end{eqnarray}
in order to prove that:
\begin{eqnarray}
\label{goal1}
|I_n-L_n|\leq |I_n-J_n|+|J_n-K_n|+|K_n-L_n|\underset{n\to\infty}{\to}0.
\end{eqnarray}

\noindent{\bf (i) Proof of $\underset{n\to\infty}{\lim}|I_n-J_n|=0$:} As $\overline{F}_{\bf t}$ is a bounded function, there exists $c>0$ such that $\frac{W_n^{(\alpha),\overline{F}_{\bf t}}}{D_n^{(\alpha)}}\leq c\frac{W_n^{(\alpha)}}{D_n^{(\alpha)}} $, then
\begin{eqnarray}
\nonumber |I_n-J_n|&\leq& c ||F||n\E_{\Q^{(\alpha)}}\Big(\frac{\1_{E_n^c}+\1_{A_n^c}} {h_\alpha(V(w_n^{(\alpha)}))}\frac{W_n^{(\alpha)}}{D_n^{(\alpha)}}\Big)
\\
\label{combin2} &\leq &  c' n \E_{\Q^{(\alpha)}}\Big(\left(\frac{\1_{E_n^c}+\1_{A_n^c}} {(h_\alpha(V(w_n^{(\alpha)}))}\right)^2 \Big)^\frac{1}{2}\E_{\Q^{(\alpha)}}\Big(\left(\frac{W_n^{(\alpha)}}{D_n^{(\alpha)}}\right)^2 \Big)^\frac{1}{2}.
\end{eqnarray}
From Lemma 4.4 \cite{AShi11} (pp15) (recall that $\Q^{(\alpha)}(A_n^c)+\Q^{(\alpha)}(E_n^c)$ converges to zero) and Lemma 4.3 in \cite{AShi11} (pp14) we can affirm that
\begin{eqnarray}
\label{combin1} \E_{\Q^{(\alpha)}}\Big(\left(\frac{\1_{E_n^c}+\1_{A_n^c}} {h_\alpha(V(w_n^{(\alpha)}))}\right)^2\Big)^\frac{1}{2}=o(\frac{1}{\sqrt{n}}),\qquad \E_{\Q^{(\alpha)}}\Big(\left(\frac{W_n^{(\alpha)}}{D_n^{(\alpha)}}\right)^2\Big)^\frac{1}{2} =O(\frac{1}{\sqrt{n}}).
\end{eqnarray}
Combining (\ref{combin1}) with (\ref{combin2}) we get (i).
\hfill $\Box$
\\

\noindent{\bf (ii) Proof of $\underset{n\to\infty}{\lim}|J_n-K_n|=0$:}  By the triangular inequality observe that $|J_n-K_n|$ is smaller
\begin{eqnarray*}
n \E_{\Q^{(\alpha)}}\Big( \frac{ \1_{ E_n}}{h_\alpha(V(w_n^{(\alpha)}))}\frac{|W_n^{(\alpha),\overline{F}_{\bf t},[k_n,n]}|}{D_n^{(\alpha)}}\Big)  +  n\E_{\Q^{(\alpha)}}\Big( \frac{ \1_{E_n}}{h_\alpha(V(w_n^{(\alpha)}))}\left|\frac{ W_n^{(\alpha),\overline{F}_{\bf t},[k_n,n]}}{D_n^{(\alpha)}} -\frac{ W_n^{(\alpha),\overline{F}_{\bf t},[k_n,n]}}{D_n^{(\alpha),[0,k_n)}}\right|   \Big)
\end{eqnarray*}
Recalling (\ref{4.3}), on $E_n$ we have $| W_n^{(\alpha),\overline{F}_{\bf t},[k_n,n]}|\leq cD_n^{(\alpha),[k_n,n]}\leq \frac{c}{n^2}$, then we deduce that
\begin{eqnarray*}
\frac{| W_n^{(\alpha),\overline{F}_{\bf t},[k_n,n]}|}{D_n^{(\alpha)} }\leq \frac{c}{D_n^{(\alpha)}n^2},\quad \text{and }\,\,  \frac{| W_n^{(\alpha),\overline{F}_{\bf t},[k_n,n]}|}{ D_n^{(\alpha),[0,k_n)}} \frac{ D_n^{(\alpha),[k_n,n]}}{D_n^{(\alpha)}}   \leq \frac{c^2}{D_n^{(\alpha)}n^2} ,\qquad \text{on } E_n.
\end{eqnarray*}
Finally it stems that $|J_n-K_n|$ is smaller
\begin{eqnarray*}
\frac{1}{n  }\E_{\Q^{(\alpha)}}\Big( \frac{\1_{E_n} }{  h_\alpha(V(w_n^{(\alpha)})) D_n^{(\alpha)}}   \Big) \leq c\frac{1}{nk_n^\frac{1}{6}} \E_{\Q^{(\alpha)}}\Big( \frac{1 }{   D_n^{(\alpha)}}\Big)=c\frac{1}{nk_n^\frac{1}{6}}=o(\frac{1}{n}),
\end{eqnarray*}
which concludes the proof of (ii).\hfill$\Box$
\\

\noindent{\bf (iii) Proof of $\underset{n\to\infty}{\lim}|K_n-L_n|=0$:} Recall (\ref{Diane2}) and (\ref{Diane4}), first observe that  \begin{eqnarray*}
 \1_{\{V(w_{k_n}^{(\alpha)})\in[k_n^{\frac{1}{3}},k_n]\}}- \1_{E_n}=\1_{E_n^c}\1_{\{V(w_{k_n}^{(\alpha)})\in[k_n^{\frac{1}{3}},k_n]\}}.
 \end{eqnarray*} 
Then let us introduce the $\sigma$-field
\begin{equation}
\mathcal{G}_p:=\sigma\left((V(w_k^{(\alpha)}))_{k\leq p},V(u)\,  \text{for }u\in \mathbb{T}\text{ such that }\exists z\in \underset{k\in [1, p]}{\bigcup}R_k^{(\alpha)} \text{ and } u\geq z \right),\quad p\in \N.
\end{equation}
Clearly $W_n^{(\alpha),[0,k_n)}$ and $D_n^{(\alpha),[0,k_n)} $ are measurable with respect to $\mathcal{G}_{k_n}$, thus 
\begin{eqnarray}
\nonumber |K_n-L_n| & \leq &          cn\E_{\Q^{(\alpha)}}\Big(\frac{W_n^{(\alpha),[0,k_n)}}{D_n^{(\alpha),[0,k_n)}}\frac{\1_{\{V(w_{k_n}^{(\alpha)})\in[k_n^{\frac{1}{3}},k_n]\}}\1_{E_n^c}}{h_\alpha(V(w_n^{(\alpha)}))}\Big)
\\
\label{THomas1} &=& cn\E_{\Q^{(\alpha)}} \Big(\frac{W_n^{(\alpha),[0,k_n)}}{D_n^{(\alpha),[0,k_n)}}\1_{\{V(w_{k_n}^{(\alpha)})\in[k_n^{\frac{1}{3}},k_n]\}}    \E_{\Q^{(\alpha)}} \Big[\frac{\1_{E_n^c}}{h_\alpha(V(w_n^{(\alpha)}))}\big|\mathcal{G}_{k_n} \Big]\Big).
\end{eqnarray}
Moreover, by the branching property, conditionally at $V(w_{k_n}^{(\alpha)})$, $\frac{\1_{E_n^c}}{h_\alpha(V(w_{n}^{(\alpha)}))}$ is independent of $\mathcal{G}_{k_n}$, so the term in (\ref{THomas1}) is equal to
\begin{eqnarray*}
&&cn\E_{\Q^{(\alpha)}} \Big(\frac{W_n^{(\alpha),[0,k_n)}}{D_n^{(\alpha),[0,k_n)}}\1_{\{V(w_{k_n}^{(\alpha)})\in[k_n^{\frac{1}{3}},k_n]\}}    \E_{\Q^{(\alpha)}} \Big[\frac{\1_{E_n^c}}{h_\alpha(V(w_n^{(\alpha)}))}\big|\sigma(V(w_{k_n}^{(\alpha)})) \Big]\Big)
\\
&\leq & cn\E_{\Q^{(\alpha)}} \Big(\frac{W_n^{(\alpha),[0,k_n)}}{D_n^{(\alpha),[0,k_n)}}\1_{\{V(w_{k_n}^{(\alpha)})\in[k_n^{\frac{1}{3}},k_n]\}} \Big)\underset{u\in[k_n^{\frac{1}{3}},k_n]}{\sup}\E_{\Q^{(\alpha)}} \Big( \frac{ \1_{E_n^c}}{h_\alpha(V(w_{n }^{(\alpha)}))}\big| V(w_{k_n}^{(\alpha)})=u \Big).
\end{eqnarray*}
According to \cite{AShi11} (see (4.9) p21) we know that
\begin{equation}
\label{3.28}
\underset{n\to\infty}{\limsup}\,\, \sqrt{n} \E_{\Q^{(\alpha)}} \Big(\frac{W_n^{(\alpha),[0,k_n)}}{D_n^{(\alpha),[0,k_n)}}\1_{\{V(w_{k_n}^{(\alpha)})\in [k_n^{\frac{1}{3}},k_n]\}} \Big)  \leq\theta.
\end{equation}
Furthermore by the Cauchy-Schwartz inequality, for any $u\in [k_n^{\frac{1}{3}},k_n]$,
\begin{eqnarray*}
\E_{\Q^{(\alpha)}} \Big( \frac{ \1_{E_n^c} }{h_\alpha(V(w_{n}^{(\alpha)}))}\big| V(w_{k_n}^{(\alpha)})=u  \Big)\leq \E_{\Q^{(\alpha)}}  \Big( \frac{1}{[h_\alpha(V(w_{n}^{(\alpha)}))]^2} \big| V(w_{k_n}^{(\alpha)})=u  \Big)^\frac{1}{2}  \Q^{(\alpha)}(E_n^c \big| V(w_{k_n}^{(\alpha)})=u)^{\frac{1}{2}}.
\end{eqnarray*}
From (\ref{eq9.1}), we have $\underset{u\in[k_n^{\frac{1}{3}},k_n]}{\sup} \Q^{(\alpha)}(E_n^c \big| V(w_{k_n}^{(\alpha)})=u)\to 0$ when $n$ goes to infinity. Concerning the first term, according to (\ref{eq3.3}), for any $u\in [k_n^\frac{1}{3},k_n]$ we have:
\begin{eqnarray*}
&&\E_{\Q^{(\alpha)}}  \Big( \frac{1}{[h_\alpha(V(w_{n}^{(\alpha)}))]^2} \big| V(w_{k_n}^{(\alpha)})=u  \Big)= \E_{\Q^{(\alpha+u)}} \Big( \frac{1}{[h_{\alpha +u}(V(w_{n-k_n}^{(\alpha+u)}))]^2}  \Big)
\\
&&=\frac{1}{h_{\alpha+u}(0)}\E \Big( \frac{  \1_{\{ \underline{S}_{n-k_n}\geq -(\alpha+u)\}}}{h_{\alpha+u}(S_{n-k_n})}  \Big)
\\
&&\leq   \overset{\sqrt{n}}{\underset{k=0}{\sum}}\frac{c(k+1)^{-1}}{h_{\alpha+u}(0)} \P_{\alpha+u} \Big( \underline{S}_{n-k_n}\geq 0,\, S_{n-k_n}\in [k,k+1]  \Big)+  \frac{cn^{-\frac{1}{2}}}{h_{\alpha+u}(0)}\P_{\alpha+u} \Big( \underline{S}_{n-k_n}\geq 0  \Big).
\end{eqnarray*}
Recalling that $k_n=o(n^{\frac{1}{2}})$ and using (\ref{lldda}) and (\ref{Aine}), we get that for any $n\in \N$ large enough and any $u\in [k_n^{\frac{1}{3}},k_n]$, 
\begin{equation}
\label{3.28comp}\E_{\Q^{(\alpha)}}  \Big( \frac{1}{[h_\alpha(V(w_{n}^{(\alpha)}))]^2} \big| V(w_{k_n}^{(\alpha)})=u  \Big)\leq c \frac{h_\alpha(u)}{h_\alpha(u)}( n^{-\frac{3}{2}}\overset{\sqrt{n}}{\underset{k=0}{\sum}}1+ n^{-1}) \leq \frac{c'}{n}.
\end{equation}
Finally by combining (\ref{3.28comp}) and (\ref{3.28}) we obtain that $\underset{n\to\infty}{\lim}|K_n-L_n|=0$.\hfill$\Box$
\\

It remains to prove that  $L_n\to 0$ when $n$ goes to infinity. By using the Markov property (assuming $n$ large enough such that $k_n\leq \underset{i\in[1,d]}{\min}t_i\, n$), we get that
\begin{eqnarray}
\nonumber L_n & =&     n\E_{\Q^{(\alpha)}} \Big(\frac{W_n^{(\alpha),\overline{F}_{\bf t},[0,k_n)}}{D_n^{(\alpha),[0,k_n)}}\1_{\{V(w_{k_n}^{(\alpha)})\in [k_n^{\frac{1}{3}},k_n]\}}\E_{\Q^{(\alpha)}} \Big[ \frac{\overline{F}_{\bf t}(V_{\bf t}(w_n^{(\alpha)}))} {h_\alpha(V(w_n^{(\alpha)}))} \big| \mathcal{G}_{k_n}   \Big]\Big)
\\
\nonumber &=&  n\E_{\Q^{(\alpha)}} \Big(\frac{W_n^{(\alpha),\overline{F}_{\bf t},[0,k_n)}}{D_n^{(\alpha),[0,k_n)}}\1_{\{V(w_{k_n}^{(\alpha)})\in [k_n^{\frac{1}{3}},k_n]\}}\E_{\Q^{(\alpha)}} \Big[ \frac{\overline{F}_{\bf t}(V_{\bf t}(w_n^{(\alpha)}))} {h_\alpha(V(w_n^{(\alpha)}))} \big|V(w_{k_n}^{(\alpha)}) \Big] \Big)
\\
\nonumber&\leq & cn\E_{\Q^{(\alpha)}} \Big(\frac{W_n^{(\alpha),[0,k_n)}}{D_n^{(\alpha),[0,k_n)}}\1_{\{V(w_{k_n}^{(\alpha)})\in[k_n^{\frac{1}{3}},k_n]\}} \Big)\times \underset{u\in[k_n^{\frac{1}{3}},k_n]}{\sup}\left|\E_{\Q^{(\alpha)}} \Big( \frac{\overline{F}_{\bf t}(V_{\bf t}(w_n^{(\alpha)}))} {h_\alpha(V(w_n^{(\alpha)}))} \big|V(w_{k_n}^{(\alpha)})=u  \Big)\right|
\\
\label{thosh} && \qquad  \leq c'\sqrt{n} \underset{u\in[k_n^{\frac{1}{3}},k_n]}{\sup}\left|\E_{\Q^{(\alpha)}} \Big( \frac{\overline{F}_{\bf t}(V_{\bf t}(w_n^{(\alpha)}))} {h_\alpha(V(w_n^{(\alpha)}))} \big|V(w_{k_n}^{(\alpha)})=u  \Big)\right|,
\end{eqnarray}
where we have used (\ref{3.28}) in the last inequality. Recalling the definition of $V_{\bf t}$ in (\ref{1.10}) using the Markov property at time $k_n$, then (\ref{eq3.3}) we can affirm that  for any $u\in [k_n^{\frac{1}{3}},k_n]$, the expectation in (\ref{thosh}) is equal to 
\begin{eqnarray*}
&&\E_{\Q^{(u+\alpha)}}\left( \frac{\overline{F}_{\bf t}\left( \frac{1}{\sqrt{n}}[V(w_{\lfloor nt_i\rfloor-k_n}^{(u+\alpha)})+ (nt-\lfloor nt \rfloor) (V(w_{\lfloor nt_i\rfloor-k_n+1}^{(u+\alpha)}) -V(w_{\lfloor nt_i\rfloor-k_n}^{(u+\alpha)})   )]_{i\in [1,d]} \right)}{h_{\alpha +u}(V(w_{n-k_n}^{(\alpha+u)}))}   \right)
\\
&&=\frac{1}{h_{\alpha+u}(0)}   \E\left( \1_{\{ \underline{S}_{n-k_n}\geq -(u+\alpha)\}} \overline{F}_{\bf t}\left( [\frac{1}{\sqrt{n}}S_{\lfloor nt_i\rfloor-k_n} + (nt-\lfloor nt \rfloor) (S_{\lfloor nt_i\rfloor-k_n+1}  -S_{\lfloor nt_i\rfloor-k_n}   )]_{i\in [1,d]} \right)\right),
\end{eqnarray*}
which we rewrite (according to (\ref{lldda})), 
\begin{equation}
\frac{1}{h_{\alpha+u}(0)}\E\left(F\left( \Delta^{(n)} {\bf S}_{\bf t}+  S_{\bf t}^{(n-k_n)} \right)\1_{\{\underline{S}_{n-k_n}+u\geq -\alpha\}}\right) - \frac{\theta}{\sqrt{n}}\E(F(\sigma R_{\bf t})) +o(\frac{1}{\sqrt{n}}),
\end{equation}
with
\begin{eqnarray}
\label{3.32} S_{\bf t}^{(n-k_n)}&:= & \frac{(S_{\lfloor (n-k_n)t_1\rfloor })_{i\in[1,d]} }{\sqrt{n-k_n}}, \quad \text{and }
\\
\label{3.33} \Delta^{(n)}{\bf S}_{\bf t}&:=&  \frac{1}{\sqrt{n}}[S_{\lfloor nt_i\rfloor-k_n} + (nt-\lfloor nt \rfloor) (S_{\lfloor nt_i\rfloor-k_n+1}  -S_{\lfloor nt_i\rfloor-k_n}   )]_{i\in [1,d]}     - S_{\bf t}^{(n-k_n)}.
\end{eqnarray}
As a straight-forward consequence of the Lemma \ref{jaffuelelargi}, in the Appendix, we have that: {\it uniformly in $u\in [k_n^\frac{1}{3},k_n]$, as $n\to \infty$,
\begin{equation}
\label{finish}
\frac{1}{ h_{\alpha+u}(0) }\E\left(F\left( \Delta^{(n)} {\bf S}_{\bf t}+  S_{\bf t}^{(n-k_n)} \right)\1_{\{\underline{S}_{n-k_n}+u\geq -\alpha\}}\right)= \frac{\theta}{\sqrt{n}}  \E(F(\sigma R_{\bf t})) +o(\frac{1}{\sqrt{n}}).
\end{equation}}
Then by combining (\ref{finish}), (\ref{thosh}) and assertion (i), (ii), (iii) we obtain (\ref{L2sousQ}) 
\hfill $\Box$
\\

We turn now to the

\noindent{\it Proof of (\ref{1particule}).} Let $\epsilon>0$. From Theorem 1.1 \cite{AShi10}, we know that
\begin{equation}
\underset{|x|=n}{\inf}V(x)\to \infty,\qquad \P^* \text{ a.s}.
\end{equation} 
Then let $k=k(\epsilon)>0$ such that $\P^{*}(\Omega_k)\geq 1-\epsilon$ with $\Omega_k:= \{ \underset{|x|\geq 0}{\inf}V(x)\geq -k\}$.

From (\ref{lldda}) there exists $M=M(\epsilon)>0$ such that 
\begin{equation}
 c_0(1-\epsilon)u\leq h_0(u)\leq c_0(1+\epsilon)u,\qquad \forall u\geq M.
\end{equation}
Now we fix $\alpha=\alpha(\epsilon):= k + M$. Since $h_\alpha(u)=h_0(u+\alpha)$, we have for all vertices $x$
$$ 0<c_0(1-\epsilon)(V(x)+\alpha)\leq h_{\alpha}(V(x))\leq c_0(1+\epsilon)(V(x)+\alpha),\quad \text{on }\Omega_k. $$
We deduce that on  $\Omega_{k}$, for any $n\in \N^*$,
\begin{eqnarray}
\label{encd1}& W_n^{(\alpha)}=W_n,\quad W_n^{(\alpha),\overline{F}_{\bf t}}= \sum_{|u|=n }\ee^{-V(u)}\overline{F}_{\bf t}(V_{\bf t}(u))&\text{ and }
\\
\label{encD}& 0<c_0(1-\epsilon)(D_n+\alpha W_n) \leq D_n^{(\alpha)}\leq c_0(1+\epsilon)(D_n+\alpha W_n).&
\end{eqnarray}
%
%
%
Furthermore $\P^*$ a.s, $D_n\underset{n\to\infty}{\to} D_\infty>0$, then let $\eta(\epsilon),\, N(\epsilon)$ such that for any $n\geq N$, $\P^*(O_n)\geq 1-\epsilon$ with $O_n:= \left\{  D_n \geq \eta \right\}$.\\

Gathering all these facts we finally deduce that for any $n\geq N(\epsilon)$, 
\begin{eqnarray*}
&& \E_{\P^*}\Big(\frac{1}{W_n^2}  \Big(\underset{|z|=n}{\sum}\ee^{-V(u)} \overline{F}_{\bf t}(V_{\bf t}(z))\Big)^2 \Big) 
\\
&&\leq c\E\Big(\frac{1}{(W^{(\alpha)}_n)^2}  \Big(\underset{|z|=n}{\sum} \1_{\{\underline{V}(u)\geq -\alpha\}}\ee^{-V(u)}  \overline{F}_{\bf t}(V_{\bf t}(z)) \Big)^2\1_{\{\Omega_k\cap O_n\}} \Big)+\P^*(O_n^c)+\P^*(\Omega_k^c)
\\
&&\leq  c'\alpha\E_{\Q^{(\alpha)}} \Big(\frac{1}{D_n^{(\alpha)}(W^{(\alpha)}_n)^2}  \Big(\underset{|z|=n}{\sum}\ee^{-V(u)} \1_{\{\underline{V}(u)\geq -\alpha\}} \overline{F}_{\bf t}(V_{\bf t}(z)) \Big)^2\1_{\{ O_n\}} \Big)+2\epsilon
\\
&&\leq \frac{c''\alpha}{\eta c_0 (1-\epsilon)}\E_{\Q^{(\alpha)}} \Big( \frac{1}{(W^{(\alpha)}_n)^2} \Big(\underset{|z|=n}{\sum}\ee^{-V(u)} \1_{\{\underline{V}(u)\geq -\alpha\}} \overline{F}_{\bf t}(V_{\bf t}(z)) \Big)^2 \Big)+2\epsilon,
\end{eqnarray*}
which is smaller that $3\epsilon$ for $n$ large enough according to (\ref{L2sousQ}). This last inequality ends the proof of (\ref{1particule}).
\hfill$\Box$

\subsection{An extension of Theorem \ref{ConvL2}}
As we will see in the next section, to prove Theorem \ref{mainresult} we will need a slightly extension of Theorem \ref{ConvL2}. Formally, it corresponds to the case where $F_C(f)=\ee^{Cf(1)},\,  f\in\mathcal{C},\, C>0$:

\begin{proposition}
\label{probaexp}
Under (\ref{criticalcondition1}), (\ref{criticalcondition2}) and (\ref{extra1}), for any $C>0$, the following equality is true in $\P^*$ probability,
\begin{equation}
\label{eqprobaexp} \underset{n\to\infty}{\lim}\frac{1}{W_n}\underset{|z|=n}{\sum}\ee^{-V(z)}\ee^{C\frac{V(z)}{\sqrt{n}}}= \E(\ee^{\sigma C R_1}),
\end{equation}
where $R_1$ denotes a Brownian meander at time $1$.
\end{proposition}
{\it Proof of Proposition \ref{probaexp}.} From Theorem \ref{ConvL2} we can affirm that
\begin{equation}\underset{p\to\infty}{\lim}\underset{n\to\infty}{\lim}\frac{1}{W_n}\underset{|u|=n}{\sum}\ee^{-V(u)}  \ee^{C\frac{V(u)}{\sqrt{n}}}\1_{\{\frac{V(u)}{\sqrt{n}}\leq p\}}=\E\left(\ee^{\sigma C R_1}\right),\quad   \text{ in } \P^* \text{ probability}.
\end{equation}
So in order to prove (\ref{eqprobaexp}) it remains to show that: {\it for any $\epsilon>0$ as $n\to \infty$ then $p\to \infty$,
\begin{equation}
\label{inter}
\P^*_{(\ref{inter})}(n,p,\epsilon):= \P^*\Big( \frac{1}{W_n}\underset{|u|=n}{\sum}\ee^{-V(u)}  \ee^{C\frac{V(u)}{\sqrt{n}}}\1_{\{\frac{V(u)}{\sqrt{n}}\geq p\}}\geq \epsilon \Big) \to 0.
\end{equation}}
Let $\epsilon>0$. We choose $k=k(\epsilon)$, $\alpha=\alpha(\epsilon)$ as in the proof of (\ref{1particule}). From (\ref{eqthm1.1}) we recall that $\underset{A\to\infty}{\lim}\underset{n\in\N}{\sup}\,\, \P^*\left(U_{n,A}^c\right)=0$ with $U_{n,A}:=\{\frac{1}{W_n}\leq A\sqrt{n} \}$. Recalling also the definition of $\Omega_k$ we deduce that for $A$ and $n$ large enough, we have
\begin{eqnarray*}
\P^*_{(\ref{inter})}(n,p,\epsilon)&\leq& c\P\Big( A\sqrt{n} \underset{|u|=n}{\sum}\ee^{-V(u)}  \ee^{C\frac{V(u)}{\sqrt{n}}}\1_{\{\underline{V}(u)\geq -\alpha,\, \frac{V(u)}{\sqrt{n}}\geq p\}}\geq  \epsilon,\, \Omega_k,\, U_{n,A} \Big)+ \P^*(\Omega_k^c)+\P^*(U_{n,A}^c)
\\
&\leq & c \frac{A\sqrt{n}}{\epsilon}\E\Big( \underset{|u|=n}{\sum}\ee^{-V(u)}  \ee^{C\frac{V(u)}{\sqrt{n}}}\1_{\{\underline{V}(u)\geq -\alpha,\, \frac{V(u)}{\sqrt{n}}\geq p\}}  \Big)+ 2\epsilon
\\
&=& \frac{c A\sqrt{n}}{\epsilon} \E\left(\ee^{C\frac{S_n}{\sqrt{n}}},\,\underline{S}_n\geq -\alpha,\, S_n\geq p\sqrt{n}\right)+2\epsilon,
\end{eqnarray*}
where in the last inequality we have used the identity (\ref{2.1}). By Lemma \ref{MajSp} (Appendix) this is smaller than $\frac{c A}{\epsilon}\ee^{-\frac{p}{4}}+2\epsilon $. Finally we conclude that
\begin{equation}
\underset{p\to\infty}{\lim}\underset{n\to\infty}{\lim} \P^*\Big( \frac{1}{W_n}\underset{|u|=n}{\sum}\ee^{-V(u)}  \ee^{C\frac{V(u)}{\sqrt{n}}}\1_{\{\frac{V(u)}{\sqrt{n}}\geq p\}}\geq \epsilon \Big)= 0,
\end{equation}
which ends the proof of Proposition \ref{probaexp}.
\hfill$\Box$

\section{Proof of Theorem \ref{mainresult}}
The proof of Theorem \ref{mainresult}  is a straightforward consequence of the followings two results: {\it Assume (\ref{criticalcondition1}), (\ref{criticalcondition2}) and (\ref{extra1}). We have:
\begin{equation}
\label{Convproba}
\underset{C\to\infty}{\lim}\underset{\beta\to1,\beta<1}{\lim}\P^*\left(\left|\frac{1}{\alpha}W_{\beta,\lfloor \frac{C}{\alpha^2}\rfloor }-2D_\infty\right|>\epsilon\right)=0,\qquad \forall \epsilon>0.
\end{equation}}
and
\begin{lemme}
\label{lemcontrol} Assume (\ref{criticalcondition1}), (\ref{criticalcondition2}) and (\ref{extra1}). We have ($\beta=1-\alpha$)
\begin{equation}
\underset{C\to \infty}{\limsup}\,  \underset{\alpha\to 0}{\limsup}\, \P^*\left( \frac{1}{\alpha}\left|W_\beta-W_{\beta,\lfloor \frac{C}{\alpha^2}\rfloor }\right|\geq \epsilon\right) =0,\qquad \forall \epsilon>0.
\end{equation}
\end{lemme}

We first prove (\ref{Convproba}).
\\

\noindent{\it Proof of (\ref{Convproba}).}  For any  $\alpha:=1-\beta>0$ (small), $C>0 $ (large), let $n=n(\alpha,C):=\lfloor \frac{C}{\alpha^2}\rfloor $. Assume (\ref{extra1}), for small $\alpha$, $\Phi(\beta)= \frac{\Phi^{''}(1)}{2}\alpha^2 +o_\alpha(\alpha^2)= \frac{\sigma^2}{2}\alpha^2+o_\alpha(\alpha^2) $, then
\begin{eqnarray*}
 W_{\beta,n}&=&\underset{|u|=n}{\sum}\ee^{-V(u)}\ee^{\alpha V(u)}\ee^{-\Phi(\beta)n} =  \ee^{-\frac{C}{2} \sigma^2+ o_\alpha(1)}\underset{|u|=n}{\sum}\ee^{-V(u)}\ee^{\sqrt{C}\frac{V(u)}{\sqrt{n}}}
\\
&=&\ee^{-\frac{C}{2} \sigma^2 + o_\alpha(1)}W_n \frac{1}{W_n}\underset{|z|=n}{\sum}\ee^{-V(z)}\ee^{\sqrt{C}\frac{V(z)}{\sqrt{n}}}.
\end{eqnarray*}
For any $C>0$, $\alpha\to 0$ implies $n\to\infty$, thus by Proposition \ref{probaexp} it stems that for any $C>0$, 
\begin{equation}
\label{touvi1} \underset{\alpha\to 0}{\lim}\, \frac{1}{W_n}\underset{|z|=n}{\sum}\ee^{-\sqrt{C} V(z)}\ee^{\frac{V(z)}{\sqrt{n}}}\to \E(\ee^{\sigma \sqrt{C} R_1}),\quad \text{ in } \P^* \text{ probability.}
\end{equation}
On the other hand, by Aïdékon and Shi \cite{AShi11} 
\begin{equation}
\label{touvi2} \sqrt{n}W_n\to \sqrt{\frac{2}{\pi\sigma^2}}D_\infty,\quad \text{ in } \P^* \text{ probability}.
\end{equation}
Combining (\ref{touvi1}) and (\ref{touvi2}) we get for any $C>0$,
\begin{equation}
\label{couvi1} \underset{\beta \to 1}{\lim} W_{\beta, \lfloor \frac{C}{\alpha^2}\rfloor }= D_\infty \ee^{-\frac{C}{2} \sigma^2 }\sqrt{\frac{2}{\pi \sigma^2}} \E\left( \ee^{\sigma \sqrt{C}R_1}\right),\quad  \text{in } \P^* \text{ probability}.
\end{equation}
Since $\E(\ee^{\sigma \sqrt{C}R_1})\sim \ee^{\sigma^2\frac{C}{2}}\sqrt{C\sigma^2}\sqrt{2\pi}$ as $C\to\infty$, we get
\begin{eqnarray}
\label{couvi2} \underset{C\to \infty}{\lim} \ee^{-\frac{C}{2}\sigma^2} \sqrt{\frac{2}{\pi \sigma^2}} \E\left( \ee^{\sigma \sqrt{C}R_1}\right)=2.
\end{eqnarray}
Finally combining (\ref{couvi1}) and (\ref{couvi2}) we obtain (\ref{Convproba}).\hfill$\Box$
\\

In order to prove the Lemma \ref{lemcontrol}, a first step consists to show the following assertion: ($\blacklozenge$) {\it There exists $c_\blacklozenge>0$ and $\alpha_0<1$ such that: $\E\left(W_\beta^{1+\frac{\alpha}{2}}\right)\leq c_\blacklozenge $ for any $0<\alpha<\alpha_0$.} 
\\

\noindent{\it Proof of ($\blacklozenge$). } We recall that under  the condition (\ref{extra2}) there exists $\epsilon_0,\, \delta_- >0$ such that $  \underset{\beta\in[1-\delta_-,1]}{\sup} \E \left( W_{\beta,1}^{1+\epsilon_0} \right)$. For any $\beta \in [1-\delta_-,1]$, let $\alpha=1-\beta$, $p=1+\frac{\alpha}{2}$. The proof of ($\blacklozenge$) is similar to this one of Lemma 3 in \cite{Big79}. Let us introduce the probability measure $\Q_\beta$ defined by
\begin{equation}
\label{defQbeta} \Q_\beta:=W_\beta.\P.
\end{equation}
We refer to \cite{Lyo97}, for the proof of the existence of this probability and the so called "spine decomposition'' of $\Q_\beta$. We shall prove that there exists $\alpha_0<1,\, c_\blacklozenge>0$ such that
\begin{equation}
\underset{\alpha\in (0,\alpha_0]}{\sup} \E(W_\beta^p)= \underset{\alpha\in (0,\alpha_0]}{\sup}\E_{\Q_\beta}(W_\beta^{p-1})\leq c_\blacklozenge.
\end{equation}
Under $\Q_\beta$, we can decompose $W_\beta$ with respect to the "spine'' $(w_n)_{n\in \N}\subset \mathbb{T}$, it leads to 
\begin{equation}
\label{Deq1}
W_\beta=\underset{|u|=1}{\sum}\ee^{-\beta V(u)-\Phi(\beta)}W_\beta(u)= \ee^{-\beta V(w_1)-\Phi(\beta)}\tilde{W}^{(1)}_\beta+\underset{v\neq w_1,\,|v|=1}{\sum}\ee^{-\beta V(v)-\Phi(\beta)}W_\beta(v),
\end{equation}
with 
\begin{equation}
W_\beta(u):=\underset{k\to\infty}{\lim}\underset{|x|=n+k,x>u}{\sum}\ee^{-\beta [V(x)-V(u)]},\qquad \tilde{W}_\beta^{(1)}:=\underset{k\to\infty}{\lim}\underset{|x|=n+k,x>w_1}{\sum}\ee^{-\beta [V(x)-V(w_1)]}.
\end{equation}
By the branching property, the random variables $\tilde{W}_\beta^{(1)}$ and $(W_\beta(u))_{|v|=1,\,w\neq w_1}$ are independent, moreover $\tilde{W}_\beta^{(1)}$ is distributed as $W_\beta$ under $\Q_\beta$ whereas for any $u$, $|u|=1,\, w_1\neq u$ $W_\beta(u)$ is distributed as $W_\beta$ under $\P$. Introducing $B_i:=\underset{v\neq w_i}{\sum}\ee^{-\beta [V(v)-V(w_i)]-\Phi(\beta)}W_\beta(v)$, $i\in\N^*$, and  iterating (\ref{Deq1}) $N$ times we get
\begin{equation}
W_\beta=\overset{N-1}{\underset{k=0}{\sum}}\ee^{-\beta V(w_k)-\Phi(\beta)k}B_{k+1}+\ee^{-\beta V(w_N)-\Phi(\beta)N}\tilde{W}_\beta^{(n)} .
\end{equation}
By convexity and observing that $\tilde{W}_\beta^{(n)}$ and  $\ee^{-\beta V(w_N)-\Phi(\beta)N} $ are independent we deduce that
\begin{equation}
\label{AnA1}\E_{\Q_\beta}(W_\beta^{p-1})\leq \E_{\Q_\beta}\Big(\big[\overset{N}{\underset{k=0}{\sum}}\ee^{-\beta V(w_k)-\Phi(\beta)k}B_{k+1}\big]^{p-1}\Big)+\ee^{-(p-1)\Phi(\beta)N}\E_{\Q_\beta}(\ee^{-\beta (p-1)V(w_N)}) \E_{\Q_\beta}(W_\beta^{p-1}).
\end{equation}
Furthermore some calculations provide $ \ee^{-(p-1)\Phi(\beta)N}\E_{\Q_\beta}(\ee^{-\beta (p-1)V(w_N)})=\ee^{[\Phi(\beta p)-p\Phi(\beta)]N}.$ As in addition $\Phi(\beta p)-p\Phi(\beta)=-\frac{3}{8}\sigma^2\alpha^2+o_\alpha(\alpha^2)$, by choosing $N=\lfloor\frac{1}{\alpha^2}\rfloor$ (and $\alpha_0$ small enough) we obtain for any $\alpha\leq \alpha_0$, 
\begin{equation}
\label{AnA2} \ee^{-(p-1)\Phi(\beta)N}\E_{\Q_\beta}(\ee^{-\beta (p-1)V(w_N)})\leq \frac{1}{2}.
\end{equation}
Combining (\ref{AnA1}) and (\ref{AnA2}) lead to
\begin{eqnarray*}
\E_{\Q_\beta}(W_\beta^{p-1})&\leq & 2\E_{\Q_\beta}\Big(\big[\overset{N}{\underset{k=0}{\sum}}\ee^{-\beta V(w_k)-\Phi(\beta)k}B_{k+1}\big]^{p-1}\Big)
\\
&\leq & cN^{p-1}\times \E_{\Q_\beta}\Big(\underset{k\leq N}{\max}\,\ee^{(-\beta V(w_k)-\Phi(\beta)k)(p-1)}\underset{k\leq N}{\max}\,B_{k+1}^{p-1}\Big).
\end{eqnarray*}
Recalling $p-1=\frac{\alpha}{2}$, we have $N^{p-1}\leq c$, moreover using the Cauchy-Schwartz inequality, we get that 
\begin{equation}
\label{eqlem0}\E_{\Q_\beta}(W_\beta^{p-1})\leq c\E_{\Q_\beta}\Big( \underset{k\leq N}{\max}\,\ee^{(-\beta V(w_k)-\Phi(\beta)k)\alpha}   \Big)^\frac{1}{2}\E_{\Q_\beta}\Big( \underset{k\leq N}{\max}\,B_{k+1} ^{\alpha}   \Big)^\frac{1}{2}.
\end{equation}
We shall bound the two terms of the product, let us start by the first. We define the random walk $\eta_k:=\alpha(-\beta V(w_k)-\Phi(\beta)k)$. For any $\alpha<1$ let $t_0(\alpha)>0$ such that 
\begin{eqnarray}
\nonumber 1=\E_{\Q_\beta}(\ee^{t_0 \eta_1})=\E_{\Q_\beta}\left(\ee^{-t_0\beta\alpha V(w_1)-t_0\alpha\Phi(\beta)}\right)&=&\ee^{\Phi(\beta(1+t_0\alpha))-\Phi(\beta)(1+t_0\alpha)}
\\
\label{del} &=&\ee^{ \frac{\sigma^2}{2}\alpha^2((t_0-1)^2-1)+o_\alpha(\alpha^2)}.
\end{eqnarray}
Then according to (\ref{del}) we can choose $\alpha_0$ small enough such that $\forall \alpha\leq \alpha_0$, $t_0(\alpha)>\frac{3}{2}$. By definition of $t_0(\alpha)$ the process $(\ee^{t_0(\alpha)\eta_k})_{k\in \N}$ is a martingale with mean $1$, so by the Doob inequality we deduce that
\begin{eqnarray*}
\P\left(\underset{k\leq N}{\max}\,\eta_k>x\right)=\P\left( \underset{k\leq N}{\max}\,  \ee^{t_0(\alpha) \eta_k} \geq \ee^{t_0(\alpha)x} \right)\leq {\ee^{-t_0 x}}, \quad x\geq 0,
\end{eqnarray*}
and thus
\begin{equation}
\label{eqme1} \E\left(\underset{k\leq N}{\max}\,\ee^{ \eta_k}\right)\leq  1+\int_0^\infty\ee^u\P\left(\underset{k\leq N}{\max}\,{\eta_k}\geq u\right)du <c<\infty.
\end{equation}

Now we need to bound $\E_{\Q_\beta}\left( \underset{k\in[1, N]}{\max}\,B_{k+1} ^{\alpha}   \right)$. Let $\kappa>\frac{1}{\epsilon_0}$. Noting that $(B_k)_{k\in \N}$ is a sequence of independent random variables identically distributed, we deduce that: 
\begin{eqnarray*}
 \E_{\Q_\beta}\left( \underset{k \in[1, N]}{\max}\,B_{k+1} ^{\alpha}   \right)= \int_0^\infty {\Q_\beta}\left( \underset{k \in[1, N]}{\max}\, B_{k+1}^{\alpha}\geq t \right)dt &=&2 +\int_2^\infty {\Q_\beta}\left( \underset{k \in[1, N]}{\max}\, B_{k+1}^{\alpha}\geq t \right)dt
\\
 &\leq& 2+N\int_2^\infty {\Q_\beta}(B_1^\frac{1}{\kappa}\geq t^{\frac{1}{\kappa \alpha}})dt.
\end{eqnarray*}
Then by a trivial change of variable, with $\alpha_0$ small enough, it stems that for any $\alpha\in [0,\alpha_0)$, 
\begin{eqnarray}
\nonumber \E_{\Q_\beta}\left( \underset{k \in[1, N]}{\max}\,B_{k+1} ^{\alpha}   \right)\leq 2+N \alpha\int_{2^{1/(\kappa\alpha)}}^\infty  {\Q_\beta}(B_1 \geq u)\frac{u^{\kappa\alpha}}{u}du &\leq & 2+\frac{2N\kappa \alpha }{2^{1/(\kappa\alpha)}}\int_{ 2^{1/(\kappa\alpha)}}^\infty {\Q_\beta}(B_1^\frac{1}{\kappa}\geq u)du
\\
\label{touitou}&\leq &2 + c\E_{\Q_\beta}(B_1^\frac{1}{\kappa}).
\end{eqnarray}
Furthermore by convexity then the branching property we have
\begin{eqnarray*}
\E_{\Q_\beta}(B_1^\frac{1}{\kappa}) = \E_{\Q_\beta}\Big(  \Big(\underset{v\neq w_1}{\sum}\ee^{-\beta V(u)-\Phi(\beta)}W_\beta^{(u)}\Big)^{\frac{1}{\kappa}}\Big)&\leq&\E_{\Q_\beta}\Big(   \E_{\Q_\beta}\Big[\underset{v\neq w_1}{\sum}\ee^{-\beta V(u)-\Phi(\beta)}W_\beta^{(u)}|\mathcal{F}_1\Big] ^{\frac{1}{t}}\Big)
\\
&=&\E_{\Q_\beta}\Big(  \Big(\underset{v\neq w_1}{\sum}\ee^{-\beta V(u)-\Phi(\beta)}\E(W_\beta^{(u)})\Big)^{\frac{1}{t}}\Big).
\end{eqnarray*}
As $\E(W_\beta^{(u)})=\E(W_\beta)=1$ and $\kappa>\frac{1}{\epsilon_0}$, we deduce that
\begin{equation}
\label{lemeq3}\E_{\Q_\beta}(B_1^\frac{1}{\kappa}) \leq \E_{\Q_\beta}\Big((W_{1,\beta})^{\frac{1}{t}}\Big)= \E\left( W_{\beta,1}^{1+\frac{1}{\kappa}}\right)\leq \underset{\beta\in[1-\delta_-,1]}{\sup}\E\left(W_{b,1}^{1+\epsilon_0}\right).
\end{equation}
Combining (\ref{touitou}) and (\ref{lemeq3}) we conclude that there exists $c>0$ such that for any $\alpha$ small enough,
\begin{equation}
\label{lemeqqq3}
\E_{\Q_\beta}\left( \underset{k \in[1, N]}{\max}\,B_{k+1} ^{\alpha}   \right)\leq c.
\end{equation}
Finally assertion ($ \blacklozenge$) follows from (\ref{eqlem0}), (\ref{eqme1}) and (\ref{lemeqqq3}).\hfill$\Box$
\\

Now we can turn to the 

\noindent{\it Proof of Lemma \ref{lemcontrol}.} In the following  $n:=\lfloor \frac{C}{\alpha^2}\rfloor$ and $p:=1+\frac{\alpha}{2} $.  In order to avoid cumbersome notation, we will assume that $\frac{C}{\alpha^2}\in \N$. The modifications needed to handle general case are minimal and straightforward, and therefore left to the reader. For any $u\in \mathbb{T}$ such that $|u|=n$, let $W_\beta^{(u)}:=\underset{k\to\infty}{\lim}\underset{|x|=n+k,x>u}{\sum}\ee^{-[V(x)-V(u)]}\overset{law}{=}W_\beta$. Note that
\begin{equation}
\frac{1}{\alpha}|W_\beta-W_{\beta,n}|=\frac{1}{\alpha}\left|\underset{|u|=n}{\sum}\ee^{-\beta V(u)-\Phi(\beta)n}(W_\beta^{(u)}-1)\right|:=\xi_n.
\end{equation}
Fix $\epsilon>0$ and set $\tilde{\xi}_n:=\E(\xi^p|\mathcal{F}_n)$. By the Markov inequality, we have
\begin{eqnarray}
\nonumber   \P^*(\xi_n\geq \epsilon) &=& c\E\left(\1_{\xi_n^p\geq\epsilon^p}\left(\1_{\tilde{\xi}_n\geq \epsilon^{1+p}}+\1_{\tilde{\xi}_n< \epsilon^{1+p}}\right)\right)
\\
\nonumber &\leq& c \P\left(\tilde{\xi}_n\geq \epsilon^{1+p}\right)+ c\E\left(\1_{\{  \tilde{\xi}_n< \epsilon^{1+p} \}} \E\left( \1_{\{ \xi_n^p\geq\epsilon^p\}} \Big| \mathcal{F}_n\right)  \right)
\\
\label{tyou2} & {\leq}& c\P\left(\tilde{\xi}_n\geq \epsilon^{1+p}\right)+ c\epsilon.
\end{eqnarray}
We shall prove that $\underset{C\to \infty}{\limsup }\, \underset{\beta \to 1}{\limsup } \,  \P\left(\tilde{\xi}_n\geq \epsilon^{1+p}\right)=0$. As $(W_\beta^{(u)}-1)_{|u|=n}$ form a sequence of independent random variables with $0$ mean, so by Petrov \cite{Pet95} ex 2.6.20 then ($\blacklozenge$) we have
\begin{eqnarray}
\nonumber \E(\xi_n^p|\mathcal{F}_n)&\leq& \frac{2}{\alpha^p}\underset{|u|=n}{\sum}\ee^{-p\beta V(u)-p\Phi(\beta)n}\E_\P\left(|W_\beta-1|^p\right)
\\
\label{thomthom} & \leq & c_\blacklozenge \frac{2}{\alpha^p}W_{p\beta,n}\ee^{n(\Phi(p\beta)-p\Phi(\beta))}
.
\end{eqnarray}
Now let us recall the three following facts:

- For $\beta\to 1$, $\Phi(\beta)= \frac{\sigma^2}{2}\alpha^2+o_\alpha(\alpha^2)$ and $p\beta= 1-\frac{\alpha}{2}+o_\alpha(\alpha)$, thus 
\begin{eqnarray}
\nonumber\Phi(p\beta)-p\Phi(\beta)=\frac{\sigma^2}{2}[(1-p\beta)^2-p\alpha^2]+o_\alpha(\alpha^2)&=& \frac{\sigma^2}{2}[\alpha^2/4-\alpha^2]+o_\alpha(\alpha^2)
\\
\label{thom1} &=&-\frac{3}{8}\alpha^2\sigma^2+o_\alpha(\alpha^2).
\end{eqnarray}

- For $\beta\to 1$,  $\alpha^p=\alpha.\alpha^{\alpha/2}\sim \alpha$.

- Let $t(\alpha):= \frac{1-p\beta}{\alpha}$ $(\underset{\alpha\to 0}{\to }\frac{1}{2} )$. Observing that $\frac{1}{\alpha} W_{p\beta,n}=  \frac{t(\alpha)}{1-p\beta}W_{p\beta,\lfloor \frac{Ct(\alpha)}{(1-p\beta)^2}\rfloor }$, by combining (\ref{couvi1}) and (\ref{couvi2}) we can affirm that
\begin{equation}
\underset{C\to\infty}{\lim}\underset{n\to\infty}{\lim}\frac{1}{\alpha} W_{p\beta,n}=D_\infty ,\quad \text{in } \P^* \text{ probability}.
\end{equation}
Thus by combining (\ref{thomthom}) and ($\blacklozenge$) to this three assertions we get that in $\P$ probability,
\begin{eqnarray}
\nonumber \underset{C\to \infty }{\limsup}\, \underset{\beta\to 1}{\limsup} \, \E_\P(\xi_n^p|\mathcal{F}_n)&\leq& cD_\infty   \underset{C\to \infty }{\limsup}\, \ee^{-\frac{3}{8}\sigma^2 C}
 \\ 
 & = &0,\quad \text{ in } \P \text{ probability.}
\end{eqnarray}
Therefore for any $\epsilon>0$, 
\begin{equation}
\label{tyou}\underset{C\to \infty}{\limsup}\,\underset{\beta \to 1}{\limsup} \,   \P(\tilde{\xi}_n\geq \epsilon)=0.
\end{equation}
With (\ref{tyou}) and (\ref{tyou2}) we obtain Lemma \ref{lemcontrol}.
\hfill$\Box$.

\section{Appendix}
Recall that $(S_n)_{n\geq 0}$ is a centred random walk with $\E(S_1^2):=\sigma^2<\infty$. For any $t\in[0,1]$, let $S_{tn}:=S_{\lfloor tn\rfloor}$. Recall also the definition of $\Delta{\bf S}^{(n)}_{\bf t} $ in (\ref{3.33}). Let $(R_s)_{s\in[0,1]}$ a Brownian meander. The following Lemma and its proof are very similar to Lemma 2.2 in \cite{JAid09}.
\begin{lemme}
\label{jaffuelelargi}
Let $(k_n)_{n\in \N}:= (\lfloor \log n\rfloor)_{n\in \N}$. For any $\t=(t_1,...,t_d)\in(0,1]^d$ and any bounded continuous function $F:[0,\infty)^d\to\r$, we have, as $n\to\infty$,
\begin{equation}
\label{deltaS}
\P\left(\sqrt{n} \Delta{\bf S}^{(n)}_{\bf t} \notin B(0,k_n),\,\underline{S}_n\geq -x \right)=o(\frac{h_0(x)}{\sqrt{n}}),
\end{equation}
and 
\begin{equation}
\label{bolgen}
\E^x\Big(F\Big(\frac{y_1+S_{t_1n}}{(n\sigma^2)^{\frac{1}{2}}},...,\frac{y_d+ S_{t_d n}}{(n\sigma^2)^{\frac{1}{2}}}\Big)\1_{\{\underline{S}_n\geq 0\}}\Big)=\frac{\theta h_0(x)}{n^{\frac{1}{2}}}\left(\E(F(R_{\bf t}))+o(1)\right),
\end{equation}
uniformly in $(x,{\bf y})\in [0,k_n]\times B(0,k_n)$.
\end{lemme}
\noindent{\it Proof of (\ref{deltaS}).}
To prove (\ref{deltaS}) we can suppose without lost of generality that $d=1$. Let $t\in (0,1]$, recalling (\ref{3.33}) observe there exists $c>0$ such that for any $n\in \N$ large enough we have 
\begin{eqnarray*}
|\Delta^{(n)}{\bf S}_t|\leq c \frac{1}{\sqrt{n}}|S_{(n-k_n)t}-S_{n t -k_n}|+\frac{k_n}{n^\frac{3}{2}}|S_{(n-k_n)t}|.
\end{eqnarray*}
It is clear that 
\begin{equation}
\label{deltaS2}\P\left(|S_{(n-k_n)t}|\geq n ,\,\underline{S}_n\geq -x \right)\leq \P(|S_{(n-k_n)t}|\geq n )=o(\frac{1}{\sqrt{n}}).
\end{equation}
Thus we only need to prove that uniformly in $x\in [0,k_n]$, 
\begin{eqnarray}
\label{deltaS1}
\P\left(|S_{(n-k_n)t}-S_{n t -k_n}|\geq k_n ,\,\underline{S}_n\geq -x \right)&=&o(\frac{h_0(x)}{\sqrt{n}}).
\end{eqnarray}
 According to the Markov property at time $nt-k_n$ we have
\begin{eqnarray*}
\P\left(|S_{(n-k_n)t}-S_{n t -k_n}|\geq k_n ,\,\underline{S}_n\geq -x \right)&\leq& \P\left(\underline{S}_n\geq -x\right) \P\left(\max( |S_{k_n(1-t)-1}|,|S_{k_n(1-t)}|)\geq k_n \right) 
\\
&\leq& \frac{h_0(x)}{\sqrt{n}} (k_n)^{-2}\E(S_{k_n(1-t)}^2)= o(\frac{h_0(x)}{\sqrt{n}}),
\end{eqnarray*}
which gives (\ref{deltaS1}). \hfill$\Box$
\\

\noindent{\it Proof of (\ref{bolgen}).}Let ${\bf t}\in (0,1]^d$ and $F:[0,\infty)^d\to \r$ a continuous function bounded by $M>0$. For any $(x,{\bf y})\in [0,b_n]\times B(0,b_n)$ we denote $\E_{(\ref{deltaS})}^x({\bf y})$ the expectation in (\ref{deltaS}). According to \cite{AJaf11} (p11) we have: {\it for any $i\in[1,d]$, $\epsilon>0$, there exists $A(\epsilon)$ large enough such that
\begin{equation}
\label{compact}
\underset{x\in[0,k_n]}{\sup}\E^x\left(\frac{S_{t_in}}{\sigma \sqrt{n}}>A(\epsilon)\right)\leq \epsilon,\quad \underset{x\in[0,k_n]}{\sup}\E^x\left(\frac{S_{t_in}}{\sigma \sqrt{n}}>A(\epsilon),\underline{S}_n\geq 0\right)\leq h_0(x)\frac{\epsilon}{\sqrt{n}}.
\end{equation}}

Thus we can suppose that $F$ is a continuous function with compact support. By approximation, we can also assume that $F$ is Lipschitz. Let $(m_n)_{n\geq 0}$ be a sequence of integers such that $\frac{n}{m_n}$ and $\frac{m_n}{k_n^2}$ go to infinity. Decomposing $\E^x_{(\ref{deltaS})}({\bf y})$ according to the time $j$ such that $\underline{S}_j=\underline{S}_n$ gives:         
\begin{eqnarray}
\nonumber \left|\E^x_{(\ref{deltaS})}({\bf y})-\underset{j=0}{\overset{m_n}{\sum}}\E^x\big(a_n((S_j+y_i)_{i\in [0,d]},k),S_j=\underline{S}_j\geq 0\big)\right|&\leq & M\E\Big(\underset{j=0}{\overset{m_n}{\sum}}\1_{\{ S_j=\underline{S}_j\geq -x,\, \underset{l\in [j,n]}{\min}(S_l-S_j)\geq 0\}}\Big)
\\
\label{thothoui} &\leq & c(n-j+1)^{-1/2} \underset{j=0}{\overset{m_n}{\sum}}\P\Big( S_j=\underline{S}_j\geq -x \Big),
\end{eqnarray}
with
\begin{equation}
\label{55.2}
a_n({\bf z},j):=\E\Big(F\left(\frac{z_1+S_{t_1n-j}}{(n\sigma^2)^{\frac{1}{2}}},...,\frac{z_d+S_{t_d n-j}}{(n\sigma^2)^{\frac{1}{2}}}\right),\underline{S}_{n-j}\geq 0\Big),\quad {\bf z}\in \r^d.
\end{equation}

The amount in (\ref{thothoui}) is negligible, indeed according to (\ref{Aine}) we have
\begin{eqnarray}
\nonumber c\underset{j=m_n+1}{\overset{n}{\sum}} (n-j+1)^{-1/2}\P^x(S_j=\underline{S}_j\geq 0) &\leq&\overset{n}{\underset{j=m_n+1}{\sum}}  (n-j+1)^{-1/2} \times c\frac{(x+1)^2}{j^{\frac{3}{2}}}
\\
\label{22.14}&\leq &c' \frac{h_0(x)}{\sqrt{n}}\frac{ k_n}{\sqrt{m_n}} .
\end{eqnarray}
Similarly, recalling that $c_1(1+x)\leq h_0(x)\leq C_1(1+x) $ and $x\leq k_n$, observe that
\begin{eqnarray}
\nonumber \underset{j=m_n+1}{\overset{n}{\sum}}\P^x\left(S_j=\underline{S}_j\geq 0\right)&\leq& c(1+x)^2\underset{j\geq m_n}{\sum}(j+1)^{-\frac{3}{2}}
\\
\label{22.15} &\leq& ch_0(x)\frac{k_n}{\sqrt{m_n}}.
\end{eqnarray}

Going back to (\ref{22.14}) let us study $ \E\Big(F\left(\frac{z_1+S_{t_1n-j}}{(n\sigma^2)^{\frac{1}{2}}},...,\frac{z_d+S_{t_d n-j}}{(n\sigma^2)^{\frac{1}{2}}}\right),\underline{S}_{n-j}\geq 0\Big),\, j\leq m_n$. Recalling (\ref{55.2}), and the definition of $A(\epsilon)$ in (\ref{compact}) , as $F$ is Lipshitz, we have for any ${\bf z}\in B(0,k_n)$,


\begin{equation}
\label{Diancrpe}  \left|a_n({\bf z},j)-\E\Big(F\left(\frac{S_{t_1 (n-j)}}{\sigma\sqrt{n-j}},...,\frac{S_{t_d (n-j)}}{\sigma\sqrt{n-j}}\right),\underline{S}_{n-j}\geq 0\Big)\right|\leq  (1)+(2),
\end{equation}
with
\begin{eqnarray*}
(1)&:=&\frac{c8A(\epsilon) \sqrt{m_n}}{\sqrt{n}}\P\left( \underline{S}_{n-j}\geq 0 \right),
\\
(2)&:=& \sum_{i=1}^d\P\left(|{\bf z}|+ |S_{t_i n-j}-S_{t_i(n-j)}|+ \frac{j}{n}S_{t_i (n-j)})\geq 8 A(\epsilon)\sqrt{m_n},\,\underline{S}_{n-j}\geq 0\right).
\end{eqnarray*}
According to (\ref{lldda}), for $n$ large enough ($m_n=o(n)$) we have $(1)\leq c'' \frac{A(\epsilon)\sqrt{m_n}}{n}$. Term (2) is quite similar to the expectation in (\ref{deltaS}). By using the Markov property and (\ref{compact}), we deduce that for any ${\bf z}\in B(0,k_n),\, j\leq m_n$, 
\begin{eqnarray*}
(2)&\leq& c \sum_{i=1}^d \big[\P\left(\underline{S}_{t_i(n-j)-j}\geq 0\right) \P\left(\max( |S_{j(1-t_i)-1}|,|S_{j(1-t_i)}|)\geq A(\epsilon)\sqrt{m_n} \right)
\\
  & &\qquad +\P\left( S_{t_i (n-j)})\geq A(\epsilon) \sqrt{n} ,\,\underline{S}_{n-j}\geq 0\right)\big]\leq \frac{c\epsilon}{\sqrt{n}}.
\end{eqnarray*}
Finally  we deduce that for $n$ large enough we have
\begin{equation}
 \label{22.16} \underset{{\bf z}\in B(0,m_n),\, j\leq m_n}{\sup}  \left|a_n({\bf z},j)-\E\Big(F\left(\frac{S_{t_1 (n-j)}}{\sigma\sqrt{n-j}},...,\frac{S_{t_d (n-j)}}{\sigma\sqrt{n-j}}\right),\underline{S}_{n-j}\geq 0\Big)\right|\leq  \frac{c\epsilon}{\sqrt{n}}.
 \end{equation}
Furthermore we know that $(\frac{S_{\lfloor nt\rfloor}}{\sigma \sqrt{n}})_{t\in[0,1]}$ conditionally to $\underline{S}_{n}\geq 0$ converges under $\P$ to the Brownian meander \cite{Bol78}. It implies that there exists $(\eta_j)_{j\geq 0}$ tending to zero such that
\begin{equation}
\label{22.17}
\left|\E\Big(F\left(\frac{S_{t_1 (n-j)}}{\sigma\sqrt{n-j}},...,\frac{S_{t_d (n-j)}}{\sigma\sqrt{n-j}}\right),\underline{S}_{n-j}\geq 0\Big)-\P(\underline{S}_{n-j}\geq 0)\E(F(R_{\t}))\right|\leq \eta_{n-j}\P(\underline{S}_{n-j}\geq 0).
\end{equation}
Let $\epsilon>0$. For $n$ large enough and $k\leq m_n$, we have from (\ref{lldda}), $|\P(\underline{S}_{n-j}\geq 0))-\frac{\theta}{\sqrt{n}}|\leq \frac{\epsilon}{\sqrt{n}}.$ Combined with (\ref{22.16}) and (\ref{22.17}), for $n$ large enough and any ${\bf z}\in B(0,k_n),\, j\leq m_n$, this gives
\begin{eqnarray*}
 &&\left|a_n({\bf z},j)  -  \P(\underline{S}_{n-j}\geq 0)\E(F(R_{\bf t}))\right|\leq c\frac{\epsilon}{\sqrt{n}}.
\end{eqnarray*}
for $n$ greater than some $n_1$, $j\leq m_n$, ${\bf z}\in B(0,k_n)$. We use this inequality for every $k=0...,m_n$ and we obtain
\begin{eqnarray*}
\underset{j=0}{\overset{m_n}{\sum}}\left| \E^x\big(a_n((S_j+y_i)_{i\in [0,d]},j) -  \frac{\theta}{\sqrt{n}}\E(F(R_{\t})),S_j=\underline{S}_j\geq 0\Big)\right| &\leq & c\frac{\epsilon}{\sqrt{n}}\underset{j=0}{\overset{m_n}{\sum}}\P^x(S_j=\underline{S}_j\geq 0)
\\
&\leq& c\frac{\epsilon}{\sqrt{n}}h_0(x).
\end{eqnarray*}
Together with (\ref{22.14}) and (\ref{22.15}), it yields that there exists $n_0\geq 0$ such that for any $n\geq n_0$, $x\in     [0,k_n]$, ${\bf y}\in B(0,k_n)$ we have
\begin{eqnarray*}
&&\left| \E^x\Big(F\Big(\frac{y_1+S_{t_1n}}{(n\sigma^2)^{\frac{1}{2}}},...,\frac{y_d+ S_{t_d n}}{(n\sigma^2)^{\frac{1}{2}}}\Big)\1_{\{\underline{S}_n\geq 0\}}\Big)-\frac{\theta h_0(x)}{n^{\frac{1}{2}}}\E(f(R_{\bf a})) \right|
\\
&&\leq  c' \frac{h_0(x)}{\sqrt{n}}\frac{ k_n}{\sqrt{m_n}}+   c\frac{\epsilon}{\sqrt{n}}h_0(x)\leq \frac{c\epsilon h_0(x)}{\sqrt{n}},
\end{eqnarray*}
which yields (\ref{bolgen}). \hfill$\Box$
\\
\\

The following lemma is a consequence of \cite{Car05} (pp 8). 
\begin{lemme}
\label{MajSp}
Assume (\ref{criticalcondition1}), (\ref{criticalcondition2}) and (\ref{extra1}). Let $(S_n)_{n\geq 0}$ be the centered random walk defined in (\ref{defSn}). For any $C,\, \alpha>0$, there exist $c(\alpha,C),\,n_0>0$ such that for any $p,n\geq n_0$, 
\begin{equation}
\label{coupiou} \E\left(\ee^{\frac{CS_n}{\sqrt{n}}};\,\underline{S}_n\geq -\alpha,\, S_n\geq p\sqrt{n}\right)\leq \frac{c}{\sqrt{n}}\ee^{-p/4}.
\end{equation}
\end{lemme}
{\it Proof of Lemma \ref{MajSp}.} Fix $C,\, \alpha>0$. According to (\ref{extra1}) we have $\Phi(1-\theta)=\E(\ee^{\theta S_1})= 1+\sigma^2\frac{\theta^2}{2}+o(\theta^2)$. We deduce that there exists $n_0=n_0(C)$ such that for any $n,\, p\geq n_0$
\begin{equation}
\label{keepmind}
\E\left(\ee^{\frac{C S_k}{\sqrt{n}}}\right)=\big[\Phi(1-\frac{C}{\sqrt{n}})\big]^n\leq c,\qquad \P\left(\frac{S_n}{\sqrt{n}}\geq p\right) \leq \frac{c}{\ee^p}.\end{equation}
 Decomposing the expectation in (\ref{coupiou}) ($:=\E{(\ref{coupiou})}$) according to the time $k$ such that $\underline{S}_k=\underline{S}_n$ yields
\begin{eqnarray}
\nonumber  \E(\ref{coupiou})&=&\sum_{k=0}^n\E\left(\1_{\{\underline{S}_k=S_k\geq -\alpha\}}\E\left( \ee^{\frac{C(x+S_{n-k})}{\sqrt{n}}}; \underline{S}_{n-k}\geq 0,\, S_{n-k}\geq p\sqrt{n}-x\right)_{\big|x=S_k}\right)
\\
\nonumber &\leq & \sum_{k=0}^n\P\left(\underline{S}_k=S_k\geq -\alpha\right) \E\left( \ee^{\frac{C S_{n-k}}{\sqrt{n}}}; \underline{S}_{n-k}\geq 0,\, S_{n-k}\geq p\sqrt{n}\right)
\\
\label{thoithoi} &\leq& (\ref{thoithoi})_1+ (\ref{thoithoi})_2,
\end{eqnarray}
with
\begin{eqnarray}
\nonumber (\ref{thoithoi})_1:= \sum_{k=\frac{n}{2}}^n\P\left(\underline{S}_k=S_k\geq -\alpha\right)\P\left(S_{n-k}\geq p \sqrt{n}\right)  & \leq& \sum_{k=\frac{n}{2}}^n \frac{c(1+\alpha^2)}{n^\frac{3}{2}}\exp(-p ) 
\\
\label{interdit0} &\leq & \frac{c'(1+\alpha)}{\sqrt{n}}\exp(-p) ,
\end{eqnarray}
(where we have used (\ref{keepmind}), the time reversal for $(S_j)_{j\leq k}$ and (\ref{Aine})) and
\begin{equation}
 \nonumber (\ref{thoithoi})_2:= \sum_{k=0}^\frac{n}{2}\P\left(\underline{S}_k=S_k\geq -\alpha\right) \E\left( \ee^{\frac{C S_{n-k}}{\sqrt{n-k}}}; \underline{S}_{n-k}\geq 0,\, S_{n-k}\geq p\sqrt{n-k}\right).\qquad\qquad\qquad \quad \,
\end{equation}
Now let us study for any $n\in \N$, 
\begin{equation}
\E\left( \ee^{\frac{C S_{n}}{\sqrt{n}}}; \underline{S}_{n}\geq 0,\, S_{n}\geq p\sqrt{n}\right).
\end{equation}
Following  Caravenna \cite{Car05} (pp 5), we define $(T_k,H_k)$ the strict ascending ladder variables process associated to the random walk $(S_n)_{n\in \N}$. Then according to (3.1) in \cite{Car05} we have:
\begin{eqnarray*}
\P\left(\frac{S_n}{\sqrt{n}}\in dx ;\,   \underline{S}_{n}\geq 0\right)&=&\frac{1}{n}\overset{n-1}{\underset{m=0}{\sum}}\int_{[0,\sqrt{n}x)}\left(\overset{\infty}{\underset{k=0}{\sum}}\P\left(T_k=m,\,H_k\in dz\right)\right)\P\left(S_{n-m}\in \sqrt{n}dx-z\right)
\\
&=&\frac{\sqrt{n}}{n}\int_{[0,1)\times[0,x)}d\nu_n(\alpha,\beta)\P\left(\frac{S_{n(1-\alpha)}}{\sqrt{n}}\in dx-\beta\right),
\end{eqnarray*}
where $\nu_n$ is the finite measure on $[0,1)\times[0,\infty)$ defined by $\nu_n(A):=\frac{1}{\sqrt{n}}\overset{\infty}{\underset{k=0}{\sum}}\P\left(\left(\frac{T_k}{n},\frac{H_k}{\sqrt{n}}\right)\in A\right).$ Applied to our case it gives,
\begin{eqnarray*}
\E\left(\ee^{C\frac{S_n}{\sqrt{n}}}\1_{\{S_n\geq p\sqrt{n}\}};  \underline{S}_{n}\geq 0\right)&=&\frac{\sqrt{n}}{n}\int_{[p,\infty)}\int_{[0,1)\times[0,x)}d\nu_n(\alpha,\beta)\ee^{Cx}\P\left(\frac{S_{n(1-\alpha)}}{\sqrt{n}}\in dx-\beta\right)
\\
&=&n^{-\frac{1}{2}}\int_{[0,1)\times[0,\infty)}\int_{[p,\infty)} \ee^{Cx}\P\left(\frac{S_{n(1-\alpha)}}{\sqrt{n}}\in dx-\beta\right)\1_{\{\beta\leq x\} }d\nu_n(\alpha,\beta)
\\
&=& n^{-\frac{1}{2}}\int_{[0,1)\times[0,\infty)} \E\left(\ee^{C\frac{S_{n(1-\alpha)}}{\sqrt{n}}}\1_{\{\frac{S_{n(1-\alpha)}}{\sqrt{n}}\geq \max(0,p-\beta)\}} \right)d\nu_n(\alpha,\beta).
\end{eqnarray*}
Furthermore as for any $p,\, n\geq n_0$,
\begin{eqnarray*}
&&\E\left(\ee^{\frac{C S_{n(1-\alpha)}}{\sqrt{n}}}\1_{\{\frac{S_{n(1-\alpha)}}{\sqrt{n}}\geq \max(0,p-\beta)\}}\right)\leq
\\
&&\qquad\qquad\qquad  \left\{ \begin{array}{ll}   c\E\left(\ee^{ C \frac{S_n}{\sqrt{n}}}\right)\leq c',\qquad \text{if } (\alpha,\beta)\in [0,1)\times [\frac{p}{2},\infty),\quad
\\
\E\left(\ee^{\frac{ C S_{n}}{\sqrt{n}}}\1_{\{\frac{S_{n}}{\sqrt{n}}\geq \frac{p}{2}\}} \right){\leq} c'' {\ee^{-p/2}},\qquad \text{if }  (\alpha,\beta)\in[0,1)\times[0,p/2),
\end{array}\right.
\end{eqnarray*}
we deduce that
\begin{eqnarray*}
\E\left(\ee^{\frac{CS_n}{\sqrt{n}}}\1_{\{S_n\geq p\sqrt{n}\}}; \underline{S}_{n}\geq 0\right)\leq c n^{-\frac{1}{2}}\left( {\ee^{-p/2}}  \nu_n([0,1)\times[0,\infty))   +    \nu_n\left([0,1)\times[p/2,\infty)\right)\right).
\end{eqnarray*}
Let $\tau(x):=\inf\{k\geq 0,\, S_k\geq x\}$. Via some usual computations, we get that
\begin{eqnarray*}
\nu_n\left([0,1)\times[p/2,\infty)\right)&=&\frac{1}{\sqrt{n}}\overset{\infty}{\underset{k=0}{\sum}}\P\left(T_k\leq n,\, H_k\geq \frac{p}{2}\sqrt{n}\right)
\\
&=&\frac{1}{\sqrt{n}}\E\left(\1_{\{\tau(\frac{p}{2}\sqrt{n})\leq n\}}\overset{\infty}{\underset{k=0}{\sum}}\1_{\{T_k\leq n,\, H_k\geq 
\frac{p}{2}\sqrt{n}\}}\right)
\\
&\overset{Markov}{=}&\frac{1}{\sqrt{n}}\E\left(\1_{\{\tau(\frac{p}{2}\sqrt{n})\leq n\}}\E\left(\overset{\infty}{\underset{k=0}{\sum}}\1_{\{T_k\leq n-\tau(\frac{p}{2}\sqrt{n})\}}\right)\right)
\\
&\leq& \P(\tau(\frac{p}{2}\sqrt{n})\leq n)  \frac{1}{\sqrt{n}}\E\left(\overset{\infty}{\underset{k=0}{\sum}}\1_{\{T_k\leq n\}}\right)
\\
&\leq& c\ee^{-\frac{p}{4}} \nu_n([0,1)\times[0,\infty)).
\end{eqnarray*}
Finally as there exists $c_1>0$ such that $\P(T_1\geq n)\geq \frac{c_1}{\sqrt{n}},\, \forall n\geq 0 $, we deduce that $\nu_n([0,1)\times[0,\infty))\leq c$ for any $n\geq 0$ and thus
\begin{equation}
\label{interdit}  
(\ref{thoithoi})_2\leq \sum_{k=0}^{\frac{n}{2}}\P\left(\underline{S}_k=S_k\geq -\alpha\right)   \frac{c\ee^{-\frac{p}{2}}}{\sqrt{n-k} }\leq \frac{c(\alpha)}{\sqrt{n}} \ee^{-\frac{p}{2}} .
\end{equation}
Finally going back to (\ref{thoithoi}), combining (\ref{interdit0}) and (\ref{interdit}) we obtain
\begin{equation}
\E(\ref{coupiou})\leq \frac{c }{\sqrt{n}}\ee^{-\frac{p}{4}},
\end{equation}
and Lemma \ref{MajSp} follows.\hfill$\Box$


%

\paragraph*{Acknowledgement}I would like to thank Elie Aïdékon, Yueyun Hu and Olivier Zindy for introducing me the problem of the transition of $W_\beta$ and for very stimulating discussions.

\bibliographystyle{plain}
\bibliography{bibli}

\end{document}